\documentclass[lettersize,journal]{IEEEtran}

\usepackage[utf8]{inputenc}
\usepackage{graphicx}
\usepackage{hyperref}      
\usepackage{amstext,psfrag}
\usepackage[dvips]{epsfig}
\usepackage{ragged2e}  
\usepackage{gensymb}
\usepackage{graphicx}
\usepackage{epstopdf}
\usepackage{color}
\usepackage{latexsym}
\usepackage{float}
\usepackage{array}
\usepackage{tabu}
\usepackage{multirow}
\usepackage{textcomp}
\usepackage{float}
\usepackage{subfig}
\usepackage{amsmath}
\usepackage{amsthm}
\usepackage{amssymb}
\usepackage[version=4]{mhchem}
\usepackage{siunitx}
\usepackage{chngcntr}
\usepackage{apptools}
\usepackage{longtable,tabularx}
\usepackage{newcmd}

\newtheorem{Definition}{Definition}
\newtheorem{Theorem}{Theorem}
\newtheorem{Corollary}{Corollary}
\newtheorem{Proposition}{Proposition}
\newtheorem{Lemma}{Lemma}
\newtheorem*{Proof*}{Proof}

\newtheorem{Remark}{Remark}

\AtAppendix{\counterwithin{Lemma}{section}}

\usepackage{balance}
\begin{document}
	\title{Geometric Active Disturbance Rejection Control of Rotorcraft on \SE \ with Fast Finite-Time Stability}
	\author{Ningshan Wang, Reza Hamrah, Amit K. Sanyal, Mark N. Glauser$^{1}$ 
		\thanks{The authors acknowledge support from the National Science Foundation award 2132799.}
		\thanks{$^{1}$ Mechanical and Aerospace Engineering, Syracuse University, Syracuse, NY 13244, USA {\tt\small \{nwang16,rhamrah,aksanyal,mglauser\}@syr.edu}}}

	\maketitle
	
	\begin{abstract}
		This article presents a tracking control framework enhanced by an extended state observer for a rotorcraft aerial vehicle modeled as a rigid body in three-dimensional translational and rotational motions. The system is considered as an underactuated system on the tangent bundle of the six-dimensional Lie group of rigid body motions, $\SE$. The extended state observer is designed to estimate the resultant external disturbance force and disturbance torque acting on the vehicle. It guarantees stable convergence of disturbance estimation errors in finite time when the disturbances are constant and finite time convergence to a bounded neighborhood of zero errors for time-varying disturbances. This extended state observer design is based on a H\"{o}lder-continuous fast finite time stable differentiator that is similar to the super-twisting algorithm, to obtain fast convergence. A tracking control scheme that uses the estimated disturbances from extended state observer for disturbance rejection, is designed to achieve fast finite-time stable tracking control. Numerical simulations are conducted to validate the proposed extended state observer and tracking control scheme with disturbance rejection. The proposed extended state observer is compared with other existing research to show its supremacy.
	\end{abstract}
	

	\section{Introduction}\label{sec:Intro}

	Small-scale rotorcraft unmanned aerial vehicles (UAVs) have become increasingly popular in various applications, such as security and monitoring, infrastructure inspection, agriculture, wildland management, package delivery, and remote sensing. However, these UAVs are frequently exposed to dynamic uncertainties and disturbances caused by turbulence induced by airflow around structures or regions. Therefore, it is crucial to ensure robust flight control performance in such challenging environments, with guaranteed stability margins even in the presence of dynamic disturbances and uncertainties. To this end, this article describes robust tracking control schemes for a rotorcraft UAV under disturbances and uncertainties.

	Recent research articles on rotorcraft the UAV tracking control schemes use different methods to tackle the adverse effects of disturbances and uncertainties during the flight. Torrente et al.~\cite{torrente2021data} use  Gaussian processes to complement the nominal dynamics of the multi-rotor in a model predictive control (MPC) pipeline. Hanover et al.~\cite{hanover2021performance} use an explicit scheme to discretize the dynamics for the nonlinear MPC solved by optimization. Bangura et al.~\cite{bangura2017thrust} use the propeller aerodynamics as a direct feedforward term on the desired thrust to re-regulate the thrust command of the rotors. Craig et al.~\cite{craig2020geometric} implement a set of pitot tubes onto the multi-rotor aircraft to directly sense the aircraft's airspeed.  With the knowledge of propeller aerodynamic characteristics, the airspeed is then utilized to obtain the disturbance forces and torques as feedforward terms to enhance control performance. Bisheban et al.~\cite{bisheban2020geometric} use artificial neural networks to obtain disturbance forces and torques with the kinematics information of the aircraft and then use the baseline control scheme based on the work by Lee et al.~\cite{lee2010geometric} in their tracking control scheme design. The methods used in these research articles either need high computational efforts~\cite{torrente2021data, hanover2021performance,bisheban2020geometric} or require precise modeling of the aerodynamic characteristics of the rotorcraft propellers~\cite{bangura2017thrust,craig2020geometric}, to obtain satisfactory control performance against disturbances.
	
	A promising control technique to maintain the control performance against disturbances and uncertainties is active disturbance rejection control (ADRC), which can be traced back to the dissertation by Hartlieb~\cite{hartlieb1956cancellation}. In an ADRC scheme, we first obtain an estimation of the unknown disturbance from a disturbance observer (DO) or an extended state observer (ESO) and then utilize it in the control design to reject the disturbance. ADRC and ESO are formally introduced in combination in~\cite{huang2001flight}, where the ESO is used to obtain disturbance estimates for disturbance rejection. Other than ESO, disturbance observer (DO)~\cite{chen2003nonlinear}, and unknown input observer (UIO)~\cite{basile1969observability} can also give disturbance estimates for a disturbance rejection control scheme.  
	
	ADRC schemes are widely used for rotorcraft UAV control.  In the research articles by Shao et al.~\cite{shao2018robust}, the disturbance estimation from asymptotically stable (AS) ESOs is employed to enhance surface trajectory tracking control scheme for a multi-rotor UAV in the presence of parametric uncertainties and external disturbances. Liu et al.~\cite{liu2022fixed} propose fixed-time stable disturbance observers (FxTSDO) and fault-tolerance mechanisms and utilize them in their translation and attitude control scheme. Mechali et al.~\cite{mechali2021observer} present FxTS ESOs for the same purpose. Wang et al.~\cite{wang2019quadrotor} implement incremental nonlinear dynamics inversion (INDI) control combined with a sliding-mode observer (SMO) for disturbance estimation and rejection. Jia et al.~\cite{jia2022accurate} employ the disturbance model obtained by Faessler et al.~\cite{faessler2017differential}, and then estimate the drag coefficient as a parameter. This disturbance model is also employed by Moeini et al.~\cite{moeini2021exponentially}. Cui et al.~\cite{cui2021adaptive} use an adaptive super-twisting ESO for the disturbance estimation. Bhale et al.~\cite{bhale2022finite} carry out disturbance estimation with the discrete-time finite-time stable (FTS) disturbance observer by Sanyal~\cite{sanyal2022discrete}.
	
	There are several methods to ensure the stability of the above-mentioned ESO/DO designs used for rotorcraft tracking control. The linear ESO by Shao et al.~\cite{shao2018robust} is AS. Mechali et al.~\cite{mechali2021observer} use the concept of geometric homogeneity~\cite{rosier1992homogeneous} to obtain an FTS ESO. A similar method is proposed in the ESO design by Guo et al.~\cite{guo2011convergence}. The Lyapunov functions/candidates used in the ESO stability analysis by Mechali et al.~\cite{mechali2021observer} and Guo et al.~\cite{guo2011convergence} are based on Rosier~\cite{rosier1992homogeneous}, and are presented implicitly. Jia et al.~\cite{jia2022accurate}, Moeini et al.~\cite{moeini2021exponentially} and Liu et al.~\cite{liu2022fixed} use variants of the DO proposed by Chen~\cite{chen2003nonlinear}. Another method is to use the super-twisting algorithm (STA)~\cite{moreno2012strict} to obtain ESO design. Xia et al.~\cite{xia2010attitude} use this method in ESO design for spacecraft attitude control, and Cui et al.~\cite{cui2021adaptive} design an adaptive super-twisting ESO using a similar method in a multi-rotor ADRC scheme.
	
	In much of the prior literature for rotorcraft UAV attitude control with ESO/DOs for disturbance torque estimation and rejection in rotational dynamics, the attitude kinematics of the ESOs/DOs are either based on local linearization or represented using local coordinates (like Euler angles) or quaternions. Local coordinate representations can have singularity issues  (e.g., gimbal lock with Euler angles), while quaternion representations may cause instability due to unwinding~\cite{bhat2000topological,chaturvedi2011rigid}. In situations where the UAVs have to carry out aggressive maneuvers, as in rapid collision avoidance for example, disturbance estimation and rejection from such schemes may not be reliable or accurate enough for precise control of the UAV.
	
	This article presents a scheme enhanced by ESO  on \SE for rotorcraft UAVs under complex and challenging aerodynamic environments. The ESO on \SE \ estimates the disturbance forces and torques during the flight of a rotorcraft UAV in both translation and rotation. The ADRC scheme on \SE \ then incorporates the disturbance estimation from the ESO and the feedback from tracking control schemes to drive the UAV to the desired trajectory. The ESO and ADRC schemes are fast finite-time stable (FFTS), abbreviated as FFTS-ESO and FFTS-ADRC, respectively. The tracking control module with fast finite-time stability is developed based on the research article by Viswanathan et al.~\cite{viswanathan2018integrated}. The FFTS-ESO design is based on a novel H\"{o}lder-continuous fast finite-time stable differentiator (HC-FFTSD). We carry out several sets of numerical simulations to show the validity of the proposed FFTS-ESO and FFTS-ADRC schemes. 
	
	We highlight some unique contributions of this article.
	\begin{itemize}
		\item The proposed ESO is the major contribution of this article. In the proposed ESO, which is the core of the proposed ADRC scheme, the pose of the rotorcraft is represented directly on the Lie group of rigid body transformations, the special Euclidean group \SE. Unlike the ESO and DO designs reported by  Mechali et al.~\cite{mechali2021observer}, Shao et al.~\cite{shao2018robust}, and  Cui et al.~\cite{cui2021adaptive}, which use Euler angles or quaternions for attitude representation or do not include attitude kinematics, like the DO by Bhale et al.~\cite{bhale2022finite} in disturbance torque estimation, the pose of the aircraft in this article is represented in \SE to avoid kinematic singularities. We do not use local coordinates (like Euler angles) or (dual) quaternions for pose representation so that we avoid singularities due to local coordinate representations or quaternion unwinding, as reported by Bhat et al.~\cite{bhat2000topological}, and Chaturvedi et al.~\cite{chaturvedi2011rigid}.  To the best of the author's knowledge, there is no existing publication on aircraft ADRC using ESO with pose representation on $\SE$. 
		\item In the FFTS-ADRC scheme, the FFTS-ESO is based on the HC-FFTSD. The commonly used geometric homogeneity method~\cite{rosier1992homogeneous,guo2011convergence,liu2019state,wang2021holder,wang2022holder}, cannot provide a straightforward (or explicit) Lyapunov function to prove the finite-time stability of the scheme. The (implicit) form of their Lyapunov functions is by Rosier~\cite{rosier1992homogeneous}. This implicit Lyapunov function complicates the robustness analysis under measurement noise and time-varying disturbances when that analysis is essential for an ESO designed for disturbance estimation in ADRC schemes. We propose HC-FFTSD as an approach inspired by the STA~\cite{moreno2012strict,vidal2016output} of sliding-mode control (SMC). This approach gives a straightforward design of a strict Lyapunov function, which is explicit, and therefore avoids the weakness mentioned above.
		\item Based on the HC-FFTSD, the proposed FFTS-ESO schemes are both FFTS and H\textup{\"{o}}lder-continuous, unlike the common STA and other FTS schemes that use discontinuous methods like terminal sliding-mode. Therefore, the proposed FFTS-ESO avoids the potentially harmful chattering phenomenon~\cite{sanyal2015finite}, while maintaining FTS convergence. 
		\item With explicit Lyapunov function in the stability analysis, we present proof of the robustness of the proposed FFTS-ESO under time-varying disturbing forces, torques, and measurement noise. To the best of the authors' knowledge, there is no prior research on the noise robustness of ESO using Lyapunov analysis.
	\end{itemize}
	The remainder of the article is as follows. Section \ref{sec:Preliminary} presents some preliminary results that are needed to obtain sufficient conditions for the stability of the ESO and ADRC schemes. HC-FFTSD is presented, along with its stability and robustness analysis in Section \ref{sec:Differentiator}. In Section \ref{sec:Problem}, the ESO and the tracking control design problems are formulated. Section \ref{sec:ESO} describes the detailed FFTS-ESO design, which is based on the differentiator design in Section \ref{sec:Differentiator}. Section \ref{sec:Tracking} obtains the control law of FFTS-ADRC for stable tracking control on \SE \ with the estimated disturbances obtained from FFTS-ESO described in Section \ref{sec:ESO}. In Section \ref{sec:Numerical}, we present two sets of numerical simulation results to validate the proposed FFTS-ESO and FFTS-ADRC, respectively. The first set of simulations validates the proposed FFTS-ESO and compares it with other existing literature on disturbance estimation to show the supremacy of FFTS-ESO. The second set of simulations is based on the discretization provided by the Lie group variational integrator (LGVI)\cite{nordkvist2010lie} model, FFTS-ESO for disturbance estimation, and the obtained control laws of FFTS-ADRC to validate the proposed disturbance rejection control scheme. We conclude the paper, in Section \ref{sec:Conclusion}, by summarizing the results and highlighting directions for forthcoming research.

	\section{Preliminaries}\label{sec:Preliminary}
	The statements and definitions in this section are used in the technical results obtained in later sections. The statements given here give the conditions under which a continuous system is finite-time stable, fast finite-time stable, and practically finite-time stable using Lyapunov analysis, and the last statement is used in developing the main result. 
	\begin{Lemma}[Finite-time stable]\label{lem:FTS}
		\textup{\cite{bhat2000finite}} Consider the following system of differential equations,
		\begin{align}\label{eqn:Nonlinear System}
			\dot{x}(t) = f(x(t)),\ f(0)=0,\ x(0)=x_0,
		\end{align}
		where $f: \cD \rightarrow\bR^n$ is continuous on an open neighborhood $\cD \subset \bR^n$ of the origin, and let there be a continuous and differentiable function $V(x(t))$ that is positive definite. Let $V(x)$ satisfy the following inequality:
		\begin{align}\label{eqn:Preliminary FTS Lyapunov Inequality}
			\dot{V}\leq -\lambda V^{\alpha},
		\end{align}
		where $x(t)\in \cD \backslash \{0\}$, $\lambda>0$, $\alpha\in]0,1[$. Then the system \eqref{eqn:Nonlinear System} is FTS at the origin, which means $\forall x_0\in \cD$, $x$ can reach the origin in finite time. Moreover, the settling time $T$, the time needed to reach the origin, satisfies
		\begin{align}\label{eqn:Preliminary FTS Settling Time}
			T \leq \frac{V^{1-\alpha}(x_0)}{\lambda(1-\alpha)}.
		\end{align}
	\end{Lemma}
	
	\begin{Lemma}[Fast finite-time stable]\label{lem:FFTS}
		\textup{\cite{yu2005continuous}} Consider the system \eqref{eqn:Nonlinear System} and let there be a continuous and differentiable function $V(x(t))$ that is positive definite. Let $V(x)$ satisfy the following inequality:
		\begin{align}\label{eqn:Preliminary FFTS Lyapunov Inequality}
			\dot{V} \leq -\lambda_1 V-\lambda_2 V^{\alpha},
		\end{align}
		where $x(t)\in \cD \backslash \{0\}$, $\lambda_1, \lambda_2>0$, $\alpha\in]0,1[$. Then the system \eqref{eqn:Nonlinear System} is FFTS at the origin and the settling time $T$ satisfies:
		\begin{align}\label{eqn:Preliminary FFTS Settling Time}
			T \leq \frac{1}{\lambda_1(1-\alpha)}\textup{ln}\frac{\lambda_1V^{1-\alpha}(x_0)+\lambda_2}{\lambda_2}.
		\end{align}
	\end{Lemma}

	\begin{Lemma}[Practically finite-time stable]\label{lem:PFTS}
		\textup{\cite{yu2005continuous,zhu2011attitude}} Consider the system \eqref{eqn:Nonlinear System} and let there be a continuous and differentiable function $V(x)$ that is positive definite. Let $V(x)$ satisfy the following inequality:
		\begin{align}\label{eqn:Preliminary PFTS Lyapunov Inequality}
			\dot{V} \leq -\lambda_1 V-\lambda_2 V^{\alpha} +\eta,
		\end{align}
		with $x(t)\in \cD \backslash \{0\}$, $\lambda>0$, and $\alpha\in]0,1[$. Then the system \eqref{eqn:Nonlinear System} is practical finite-time stable (PFTS) at the origin, which means that
		the solution of \eqref{eqn:Nonlinear System} will converge to the following set in finite time
		\begin{align*}
			\left\{x\,  \bigg| \, V(x)\leq \textup{min}\left\{ \frac{\eta}{(1-\theta_0)\lambda_1}, \left(\frac{\eta}{(1-\theta_0)\lambda_2} \right)^{\frac{1}{\alpha}} \right\} \right\},
		\end{align*}
		where $0<\theta_0<1$. The settling time $T$ is bounded above as follows:
		\begin{align*}
			T \leq  \textup{max}&\left\{ t_0  + \frac{1}{\theta_0\lambda_1(1-\alpha)}\textup{ln}\frac{\theta_0\lambda_1V^{1-\alpha}(x_0)+\lambda_2}{\lambda_2},\right. \\
			&\left.  t_0  + \frac{1}{\lambda_1(1-\alpha)}\textup{ln}\frac{\lambda_1V^{1-\alpha}(x_0)+\theta_0\lambda_2}{\theta_0\lambda_2}\right\}.
		\end{align*}
	\end{Lemma}
	
	\begin{Lemma}\label{lem:Binom}
		\textup{\cite{hardy1952inequalities}}Let $x$ and $y$ be non-negative real numbers and let $p\in ]1,2[$. Then 
		\begin{align}\label{eqn:Bires}
			x^{\frac{1}{p}}+ y^{\frac{1}{p}} \ge (x+y)^{\frac{1}{p}}.
		\end{align}
		Moreover, the above inequality is a strict inequality if both $x$ and $y$ are non-zero. 
	\end{Lemma}
	
	\begin{Definition}\label{def:H}
		Define $H: \mathbb{R}^3 \times \mathbb{R}\rightarrow \textup{Sym(3)}$, the space of symmetric $3\times3$ matrices, as follows:
		{\em
			\begin{equation}\label{eqn:H}
				H(x,k) := I - \frac{2k}{x\Tp x} xx\Tp. 
			\end{equation}
		}
	\end{Definition}
	
	\begin{Lemma}\label{lem:Inequality Noise Robustness}
		Let $\mu\in\bR^n/\{0\}$ and $\alpha\in]0,1/2[$. Consider $\cD:\bR^n\setminus \{0,-\mu\}$ and define $\phi(x):\cD \rightarrow \bR^+$ as:
		{\em
			\begin{align}\label{eqn:Function phi}
				\begin{split}
					&\phi(x) := Y(x)\Tp Y(x), \mbox{ where } \\
					&Y(x) :={\|x\|^{-2\alpha}}x - {\|x+\mu\|^{-2\alpha}}(x+\mu).
				\end{split}
			\end{align}
		}
		The global maximum of $\phi(x)$ is at $x=-\mu/2$.
	\end{Lemma}
	We attach the proof of Lemma \ref{lem:Inequality Noise Robustness} in the appendix. 
	
	\section{H\"{o}lder-Continuous Fast Finite-Time Stable Differentiator (HC-FFTSD)}\label{sec:Differentiator}
	In this section, we design the error dynamics for the proposed ESO in Section \ref{sec:ESO} in the form of an HC-FFTSD. We analyze the stability and robustness of the proposed HC-FFTSD in this section, to support the 
	development of the ESO design in Section \ref{sec:ESO}. Theorem \ref{thr:FFTS Differentiator} covers the stability proof of the proposed HC-FFTSD.  Corollary \ref{cor:FTS Differentiator Disturbance Robustness} describes the convergence performance of the differentiator under external disturbances.  
	Corollary \ref{cor:FTS Differentiator Noise Robustness} describes the convergence performance of the differentiator under measurement noise. In the analysis that follows, $e_1 \in \bR^n$ stands for the measurement estimation error and $e_2 \in \bR^n$
	stands for the disturbance estimation error in the ESO error dynamics, respectively. In this section and the remainder of this paper, we denote the minimum and maximum eigenvalues of a matrix by $\lambda_{\min}(\cdot)$ and $\lambda_{\max}(\cdot)$, respectively.
	\begin{Theorem}\label{thr:FFTS Differentiator}
		Let $p\in]1,2[$ and $k_3>0$. Define $\phi_1(\cdot): \bR^n \rightarrow \bR^n$ and  $\phi_2(\cdot): \bR^n \rightarrow \bR^n$ as follows: 
		\begin{align}\label{eqn:phi1phi2}
			\begin{split}
				\phi_1(e_1) &=    k_3 e_1 +(e_1\Tp e_1)^\frac{1-p}{3p-2}e_1, \\
				\phi_2(e_1) &=    k_3^2 e_1 + \frac{2k_3(2p-1)}{3p-2} (e_1\Tp e_1)^\frac{1-p}{3p-2}e_1\\
				&+\frac{p}{3p-2}(e_1\Tp e_1)^\frac{2(1-p)}{3p-2}e_1. 
			\end{split}
		\end{align}
		Define the differentiator gains $k_1, k_2>0$ and  $\mathcal{A}^*\in \bR^{2\times2}$, as:
		\begin{align}\label{eqn:A}
			\mathcal{A}^* = 
			\left[
			\begin{array}{cc}
				-k_1  &  1 \\
				-k_2  & 0 \\
			\end{array}
			\right],	
		\end{align}
		which makes $\mathcal{A}^*$ a Hurwitz matrix.
		Thereafter, the differentiator design:
		\begin{align}\label{eqn:Differentiator}
			\begin{array}{ll}
				\dot{e}_1 &= -k_1 \phi_1(e_1) + e_2, \\ 
				\dot{e}_2 &= -k_2 \phi_2(e_1), 
			\end{array}
		\end{align}
		ensures that $(e_1\Tp,e_2\Tp)\in\bR^{2n}$ converges to the origin in a fast finite-time stable manner.
	\end{Theorem}
	\begin{Proof*}
		{\em
			The proof of Theorem \ref{thr:FFTS Differentiator} is based on Theorem 1 by Vida et al.~\cite{vidal2016output}, Theorem 1 by Moreno \cite{moreno2012strict} and Proposition 3 by Cruz~\cite{cruz2011uniform}. Two properties of $\phi_1$ and $ \phi_2$ are provided as follows. \\
			\textit{Property 1 (P1): The Jacobian of $\phi_1(e_1)$, denoted $\phi_1'(e_1)$, is given as follows:}
			\begin{align}\label{eqn:Property1A}
				\begin{split}
					\phi'_1(e_1)&=\frac{\textup{d}\phi_1(e_1)}{\textup{d}e_1} \\
					&= k_3 I +(e_1\Tp e_1)^\frac{1-p}{3p-2}
					\bigg[
					I-\frac{2(p-1)}{3p-2}\frac{e_1e_1\Tp}{e_1\Tp e_1}
					\bigg],
				\end{split}
			\end{align}
			\textit{so that the following identity holds:}
			\begin{align}\label{eqn:Property1B}
				\phi_2(e_1)=\phi'_1(e_1)\phi_1(e_1)
			\end{align}
			\textit{Property 2 (P2): $\phi_1'$ is a positive definite matrix, which means $\forall w\in \bR^{2n}, e_1\in \bR^n$,}
			\begin{align}\label{eqn:Property2A}
				\lambda_{\textup{min}} \{\phi’_1(e_1)\}||w||^2  \leq w\Tp\phi'_1(e_1)w\leq \lambda_{\textup{max}} \{\phi’_1(e_1)\}||w||^2. 
			\end{align}
			\textit{The maximum and minimum eigenvalues of $\phi'_1(e_1)$ employed in \eqref{eqn:Property2A} are as given below: }
			\begin{align}
				&\lambda_{\textup{max}} \{\phi'_1(e_1)\}=k_3 + (e_1\Tp e_1)^{\frac{1-p}{3p-2}}, \label{eqn:Property2B max} \\
				&\lambda_{\textup{min}} \{\phi'_1(e_1)\}=k_3 + (e_1\Tp e_1)^{\frac{1-p}{3p-2}} \frac{p}{3p-2}\label{eqn:Property2B min}. 
			\end{align}
			From Theorem 5.5 by Chen~\cite{chen1984linear}, we know that for a Hurwitz matrix $\mathcal{A}^*$ as in \eqref{eqn:A}, $\forall\, \mathcal{Q}^* \in \bR^{ 2\times 2}$ where $\mathcal{Q}^*  \succ 0$, the Lyapunov equation: 
			\begin{align}\label{eqn:Lyapunov Equation}
				(\mathcal{A}^*)\T \mathcal{P}^*+\mathcal{P}^*\mathcal{A}^* = -\mathcal{Q}^* ,
			\end{align}
			has a unique solution $\mathcal{P}^*\succ 0$. Express the positive definite matrices $\mathcal{P}^*$ and $\mathcal{Q}^*$ in components as:
			\begin{align}
				\mathcal{P}^* = 
				\left[
				\begin{array}{cc}
					p_{11}  &  p_{12} \\
					p_{12}  &  p_{22} \\
				\end{array}
				\right], \
				\mathcal{Q}^* = 
				\left[
				\begin{array}{cc}
					q_{11}  &  q_{12} \\
					q_{12}  &  q_{22} \\
				\end{array}
				\right]. \notag
			\end{align}
			With $\mathcal{P}^*$ defined as the solution to \eqref{eqn:Lyapunov Equation},  $\mathcal{A}^*$, 
			$\mathcal{P}^*$ and $\mathcal{Q}^*$ can be augmented to $\mathcal{A}, \mathcal{P}, \mathcal{Q} \in \bR^{2n\times 2n}$, as follows:
			\begin{align*}
				\begin{split}
					&\mathcal{A} = 
					\left[
					\begin{array}{cc}
						-k_1 I  &  I \\
						-k_2 I &  0 \\
					\end{array}
					\right], \\
					&\mathcal{P} = 
					\left[
					\begin{array}{cc}
						p_{11}I  &  p_{12}I \\
						p_{12}I &  p_{22}I \\
					\end{array}
					\right], 
					\mathcal{Q} = 
					\left[
					\begin{array}{cc}
						q_{11}I  &  q_{12}I \\
						q_{12}I  &  q_{22}I \\
					\end{array}
					\right]. 
				\end{split}
			\end{align*}
			The augmented matrices $\mathcal{A}, \mathcal{P}, \mathcal{Q}$ defined above also satisfy a Lyapunov equation as given below:
			\begin{align}\label{eqn:Augmented Lyapunov Equation}
				\mathcal{A}\T\mathcal{P}+\mathcal{P}\mathcal{A}=-\mathcal{Q}.
			\end{align}
			Further, the eigenvalues of $\mathcal{P}$ and $\mathcal{P}^*$, are related such that $\lambda_{\textup{min}}\{\mathcal{P}^*\} = \lambda_{\textup{min}}\{\mathcal{P}\} $, and $\lambda_{\textup{max}}\{\mathcal{P}^*\} = \lambda_{\textup{max}}\{\mathcal{P}\} $. Similar relations hold for $\mathcal{Q}$ and $\mathcal{Q}^*$.
			Thus, with $\mathcal{P}$ as the solution to \eqref{eqn:Augmented Lyapunov Equation}, we consider the following Lyapunov candidate:
			\begin{align}\label{eqn:Lyapunov FTS Differentiator}
				V(e_1,e_2)=\zeta\T \mathcal{P}\zeta, 
			\end{align}
			where $\zeta \in \bR^{2n}$ is defined as $\zeta := [\phi_1\T(e_1),e_2\T]\T$ and $\mathcal{P}$ is the augmented $\mathcal{P}^*$, which is the unique solution of \eqref{eqn:Lyapunov Equation} for a given $\mathcal{Q}^*  \succ 0$. The upper and lower bounds of the Lyapunov candidate $V$ in \eqref{eqn:Lyapunov FTS Differentiator} are as given below:
			\begin{align}\label{eqn:Lyapunov Bound 1}
				\lambda_\textup{min}\left\{\mathcal{P}\right\} \|\zeta\|^2 \leq V(e_1,e_2) \leq \lambda_\textup{max}\left\{\mathcal{P}\right\} \|\zeta\|^2.
			\end{align}
			From \eqref{eqn:Lyapunov Bound 1}, we obtain the following two inequalities:
			\begin{align}
				&\lambda_\textup{min}\left\{\mathcal{P}\right\} (e\T_1e_1)^\frac{p}{3p-2} \leq  \lambda_\textup{min}\left\{\mathcal{P}\right\} \|\zeta\|^2 \leq V(e_1,e_2), \label{eqn:Lyapunov Bound 2} \\
				&k^2_3\lambda_\textup{min}\left\{\mathcal{P}\right\} e\T_1e_1 \leq  \lambda_\textup{min}\left\{\mathcal{P}\right\}\|\zeta\|^2 \leq V(e_1,e_2). \label{eqn:Lyapunov Bound 3}
			\end{align}
			$V(e_1,e_2)$ is differentiable everywhere except the subspace $\mathcal{S}=\{(e_1,e_2)\in\bR^{2n}|e_1 = 0 \}$.
			From \eqref{eqn:Differentiator} and Property (P1), we obtain the time derivative of $\zeta$ as follows,
			\begin{align}\label{eqn:zeta Derivtive}
				\begin{split}
					\dot{\zeta} &=
					\begin{bmatrix}
						\phi'_1(e_1)\dot{e}_1  \\
						\dot{e}_2 
					\end{bmatrix} =
					\begin{bmatrix}
						\phi'_1(e_1)(-k_1\phi_1(e_1)+e_2)\\
						-k_2  \phi'_1(e_1) \phi_1(e_1)
					\end{bmatrix} \\
					& = \mathcal{D}(e_1)\mathcal{A}\zeta ,
				\end{split}
			\end{align}
			where,
			\begin{align}\label{eqn:D}
				\begin{split}
					&\mathcal{D}(e_1) = \textup{diag}[\phi'_1(e_1),\phi'_1(e_1)]\in \bR^{2n\times 2n}, \\
					&\lambda_\textup{min}\left\{\mathcal{D}(e_1)\right\} = \lambda_\textup{min}\left\{\phi'_1(e_1)\right\}.
				\end{split}
			\end{align}
			With the expression of $\dot{\zeta}$ in \eqref{eqn:zeta Derivtive}, we obtain the time derivative of $V(e_1, e_2)$ as
			\begin{align}\label{eqn:Lyapunov Derivative 1}
				\begin{split}
					\dot{V} &= \dot{\zeta}\T \mathcal{P} \zeta + \zeta\T \mathcal{P}\dot{\zeta}  \\
					&= \zeta\T ((\mathcal{D}(e_1)\mathcal{A})\T \mathcal{P}+\mathcal{P}\mathcal{D}(e_1)\mathcal{A})\zeta  \\
					&= -\zeta\T \mathcal{\bar{Q}}(e_1) \zeta. 
				\end{split}
			\end{align}
			where $\mathcal{\bar{Q}}(e_1) $ is as
			\begin{align}\label{eqn:bar Q}
				\begin{split}
					\mathcal{\bar{Q}}(e_1) &= (\mathcal{D}(e_1)\mathcal{A})\T \mathcal{P}+\mathcal{P}\mathcal{D}(e_1)\mathcal{A} \\
					&=
					\begin{bmatrix}
						\mathcal{\bar{Q}}_{11}(e_1) & \mathcal{\bar{Q}}_{12}(e_1) \\
						\mathcal{\bar{Q}}_{12}(e_1) & \mathcal{\bar{Q}}_{22}(e_1) 
					\end{bmatrix}, \\
					\mathcal{\bar{Q}}_{11}(e_1) &= 2(k_1p_{11}+k_2p_{12})\phi'_1(e_1), \\
					\mathcal{\bar{Q}}_{12}(e_1) &= (k_1p_{12}+k_2p_{22} -p_{11})\phi'_1(e_1),\\
					\mathcal{\bar{Q}}_{22}(e_1) &= -2p_{12}\phi'_1(e_1). 
				\end{split}
			\end{align}
			With \eqref{eqn:bar Q} and \eqref{eqn:Augmented Lyapunov Equation}, we obtain $\mathcal{\bar{Q}} = \mathcal{Q}\mathcal{D}(e_1) $. Afterwards, with $\mathcal{Q},\; \mathcal{D}(e_1) \succ 0$, as defined in \eqref{eqn:Augmented Lyapunov Equation} and \eqref{eqn:D}, following inequality on their eigenvalues holds:
			With $\mathcal{Q} \succ 0$ and $\mathcal{D}(e_1)\succ 0$, we obtain following inequality on their eigenvalues,
			\begin{align}\label{eqn:Eigen Inequality 1}
				\lambda_\textup{min}\left\{\mathcal{Q} \mathcal{D}(e_1)\right\} \geq \lambda_\textup{min}\left\{\mathcal{Q} \right\} \lambda_\textup{min}\left\{\mathcal{D}(e_1)\right\} >0.
			\end{align}
			With Property 2, substituting \eqref{eqn:Eigen Inequality 1} into \eqref{eqn:Lyapunov Derivative 1} we obtain
			\begin{align}\label{eqn:Lyapunov Derivative 2}
				\begin{split}
					\dot{V} &= -\zeta\T (\mathcal{Q} \mathcal{D}(e_1)) \zeta \\
					& \leq -\lambda_\textup{min}\left\{\mathcal{Q} \mathcal{D}(e_1)\right\}\zeta\T \zeta \\
					& \leq -\lambda_\textup{min}\left\{\mathcal{D}(e_1)\right\}\lambda_\textup{min}\left\{\mathcal{Q} \right\}\zeta\T \zeta 
				\end{split}
			\end{align}
			With $\lambda_\textup{min}\left\{\mathcal{D}(e_1)\right\} = \lambda_\textup{min}\left\{\phi'_1(e_1)\right\}$,  substituting \eqref{eqn:Property2B min} and \eqref{eqn:Lyapunov Bound 2} into \eqref{eqn:Lyapunov Derivative 2}, we obtain,
			\begin{align}\label{eqn:Lyapunov Derivative 3}
				\begin{split}
					\dot{V} & \leq-\Big[ k_3 + (e_1\T e_1)^{\frac{1-p}{3p-2}} \frac{p}{3p-2} \Big]\lambda_\textup{min}\left\{\mathcal{Q} \right\} \zeta\T \zeta \\
					&\leq -\frac{\lambda_\textup{min}\left\{\mathcal{Q} \right\}}{\lambda_\textup{max}\left\{\mathcal{P}\right\}}\Big[ k_3 + \Big(\frac{V}{\lambda_\textup{min}\left\{\mathcal{P}\right\}}\Big)^\frac{1-p}{p} \frac{p}{3p-2} \Big] V  \\
					& \leq -\gamma_1 V -\gamma_2 V^\frac{1}{p}, 
				\end{split}
			\end{align}
			where $\gamma_1$ and $\gamma_2$ are positive constants, defined as,
			\begin{align} \label{eqn:gamma1 gamma2}
				\begin{split}
					\gamma_1 &= k_3 \frac{\lambda_\textup{min}\left\{\mathcal{Q} \right\}}{\lambda_\textup{max}\left\{\mathcal{P}\right\}} = k_3 \frac{\lambda_\textup{min}\left\{\mathcal{Q}^* \right\}}{\lambda_\textup{max}\left\{\mathcal{P}^*\right\}};\\ 
					\gamma_2 &= \frac{\lambda_\textup{min}\left\{\mathcal{Q} \right\}\lambda_\textup{min}\left\{\mathcal{P}\right\}^\frac{p-1}{p}}{\lambda_\textup{max}\left\{\mathcal{P}\right\}} \frac{p}{3p-2} \\
					&= \frac{\lambda_\textup{min}\left\{\mathcal{Q}^* \right\}\lambda_\textup{min}\left\{\mathcal{P}^*\right\}^\frac{p-1}{p}}{\lambda_\textup{max}\left\{\mathcal{P}^*\right\}} \frac{p}{3p-2}. 
				\end{split}
			\end{align}
			Therefore, based on the inequality \eqref{eqn:Lyapunov Derivative 3}, Lemma \ref{lem:FTS} and Lemma \ref{lem:FFTS}, we conclude that the origin of the error dynamics \eqref{eqn:Differentiator} is finite-time stable and fast finite-time stable.
		}
		\qedsymbol{}
	\end{Proof*}
	\begin{Corollary}[Disturbance Robustness]\label{cor:FTS Differentiator Disturbance Robustness}
		Consider the proposed HC-FFTSD \eqref{eqn:Differentiator} in Theorem \ref{thr:FFTS Differentiator} under perturbation, $ \delta = (\delta_1\Tp, \delta_2\Tp)\Tp$, $\delta_1, \delta_2 \in \bR^n$, and $\delta$ is bounded as $||\delta||\leq \bar{\delta}$. Thereafter, the differentiator under perturbation is as
		\begin{align}\label{eqn:Differentiator Perturbation}
			\begin{array}{ll}
				\dot{e}_1 &= -k_1 \phi_1(e_1) + e_2+\delta_1, \\ 
				\dot{e}_2 &= -k_2 \phi_2(e_1) + \delta_2.
			\end{array}    
		\end{align}
		When $\gamma_1$ in \eqref{eqn:gamma1 gamma2} fulfills $\gamma_1 \geq {\lambda_\textup{max}\left\{\mathcal{P}\right\}}/ {\lambda_\textup{min}\left\{\mathcal{P}\right\}}$, \eqref{eqn:Differentiator Perturbation} is Practically Finite-Time Stable (PFTS).
	\end{Corollary}
	\begin{Proof*}
		{\em
			Consider the Lyapunov stability analysis in Theorem \ref{thr:FFTS Differentiator}. With the Lyapunov-candidate defined by \eqref{eqn:Lyapunov FTS Differentiator} and the expression of the differentiator under perturbation in \eqref{eqn:Differentiator Perturbation}, we express the time derivative of \eqref{eqn:Lyapunov FTS Differentiator} as follows:
			\begin{align}\label{eqn:Lyapunov Derivative 4}
				\dot{V} & \leq  -\gamma_1 V -\gamma_2 V^\frac{1}{p} + 2\lambda_\textup{max}\left\{\mathcal{P}\right\}\bar{\delta} ||\zeta||.
			\end{align}
			By applying Cauchy-Schwarz inequality and \eqref{eqn:Lyapunov Bound 1}, from \eqref{eqn:Lyapunov Derivative 4}, we obtain,
			\begin{align}\label{eqn:Lyapunov Derivative 5}
				\begin{split}
					\dot{V} & \leq  -\gamma_1 V -\gamma_2 V^\frac{1}{p} + 
					\lambda_\textup{max}\left\{\mathcal{P}\right\} ||\zeta||^2 + \lambda_\textup{max}\left\{\mathcal{P}\right\} \bar{\delta}^2 \\
					& \leq -\Big(\gamma_1- \frac{\lambda_\textup{max}\left\{\mathcal{P}\right\}}{\lambda_\textup{min}\left\{\mathcal{P}\right\}}\Big) V - \gamma_2 V^\frac{1}{p} + \lambda_\textup{max}\left\{\mathcal{P}\right\} \bar{\delta}^2.
				\end{split}
			\end{align}
			Therefore, according to Lemma \ref{lem:PFTS}, with inequality \eqref{eqn:Lyapunov Derivative 5}, we conclude that the system \eqref{eqn:Differentiator Perturbation}, which is the differentiator \eqref{eqn:Differentiator} under disturbance $\delta$, is practical finite time stable at the origin.
		}\qedsymbol{}
	\end{Proof*}
	\begin{Corollary}[Noise Robustness]\label{cor:FTS Differentiator Noise Robustness}
		Consider the proposed HC-FFTSD \eqref{eqn:Differentiator} in Theorem \ref{thr:FFTS Differentiator} under measurement noise $\mu$, so that $\phi_1(e_1)$ and $\phi_2(e_1)$ in \eqref{eqn:phi1phi2} are as $\phi_1(e_1+\mu)$ and $\phi_2(e_1+\mu)$. Thereafter
		\begin{align}\label{eqn:Differentiator Noise 1}
			\begin{array}{ll}
				\dot{e}_1 &= -k_1 \phi_1(e_1+\mu) + e_2 \\ 
				\dot{e}_2 &= -k_2 \phi_2(e_1+\mu), 
			\end{array}    
		\end{align}
		where $\mu$ is bounded as $||\mu||\leq \bar{\mu}$. When $\gamma_1$ in \eqref{eqn:gamma1 gamma2} fulfills $\gamma_1 \geq {\lambda_\textup{max}\left\{\mathcal{P}\right\}}/ {\lambda_\textup{min}\left\{\mathcal{P}\right\}}$, \eqref{eqn:Differentiator Noise 1} is practically finite-time stable (PFTS).
	\end{Corollary}
	\begin{Proof*}
		{\em
			From \eqref{eqn:Differentiator Noise 1}, we obtain the following expression
			\begin{align}\label{eqn:Differentiator Noise 2}
				\begin{array}{ll}
					&\dot{e}_1 = -k_1 \phi_1(e_1)+ e_2+k_1 \phi^*_1(e_1,\mu),  \\ 
					&\dot{e}_2 = -k_2 \phi_2(e_1)+k_2 \phi^*_2(e_1,\mu), \\
					&\phi^*_1(e_1,\mu) = -\phi_1(e_1+\mu) + \phi_1(e_1), \\
					&\phi^*_2(e_1,\mu) = -\phi_2(e_1+\mu) + \phi_2(e_1).
				\end{array}    
			\end{align}
			With \eqref{eqn:phi1phi2}, we obtain 
			\begin{align*}
				&\phi^*_1(e_1,\mu) = -\phi_1(e_1+\mu) + \phi_1(e_1) \\
				&= -k_3\mu - \|e_1+\mu\|^\frac{2(1-p)}{3p-2}(e_1+\mu) + \|e_1\|^\frac{2(1-p)}{3p-2}e_1 \\
				&\phi^*_2(e_1,\mu) = -\phi_2(e_1+\mu) + \phi_2(e_1)\\
				&= -k^2_3\mu - \frac{2k_3(2p-1)}{3p-2}\|e_1+\mu\|^\frac{2(1-p)}{3p-2}(e_1+\mu)\\
				& +\frac{2k_3(2p-1)}{3p-2} \|e_1\|^\frac{2(1-p)}{3p-2}e_1+\frac{p}{3p-2}\|e_1\|^\frac{4(1-p)}{3p-2}e_1\\
				&-\frac{p}{3p-2}\|e_1+\mu\|^\frac{4(1-p)}{3p-2}(e_1+\mu). 
			\end{align*}
			Therefore, according to Lemma \ref{lem:Inequality Noise Robustness}, we obtain the upper bounds of $\|\phi^*_1(e_1,\mu)\| $ and $\|\phi^*_2(e_1,\mu)\| $ as:
			\begin{align*}
				\|\phi^*_1(e_1,\mu)\|& \leq k_3 \bar{\mu} + 2^{\frac{2(p-1)}{3p-2}} (\bar{\mu})^{1-\frac{2(p-1)}{3p-2}}\\
				\|\phi^*_2(e_1,\mu)\|& \leq k^2_3 \bar{\mu} +  \frac{2k_3(2p-1)}{3p-2}2^{\frac{2(p-1)}{3p-2}}(\bar{\mu})^{1-\frac{2(p-1)}{3p-2}} \\
				&+\frac{p}{3p-2}2^{\frac{4(p-1)}{3p-2}}(\bar{\mu})^{1-\frac{4(p-1)}{3p-2}}.
			\end{align*}
			Thus, with upper bounded $\|\phi^*_1(e_1,\mu)\|$ and $\|\phi^*_1(e_1,\mu)\|$, by Corollary \ref{cor:FTS Differentiator Disturbance Robustness}, we conclude that the error dynamics \eqref{eqn:Differentiator Noise 1} is PFTS at the origin.
		}\qedsymbol{}
	\end{Proof*}
	
	\section{Problem Formulation}\label{sec:Problem}
	\subsection{Coordinate frame definition}
\begin{figure}[ht]
		\centering
		\includegraphics[width=0.6\columnwidth]{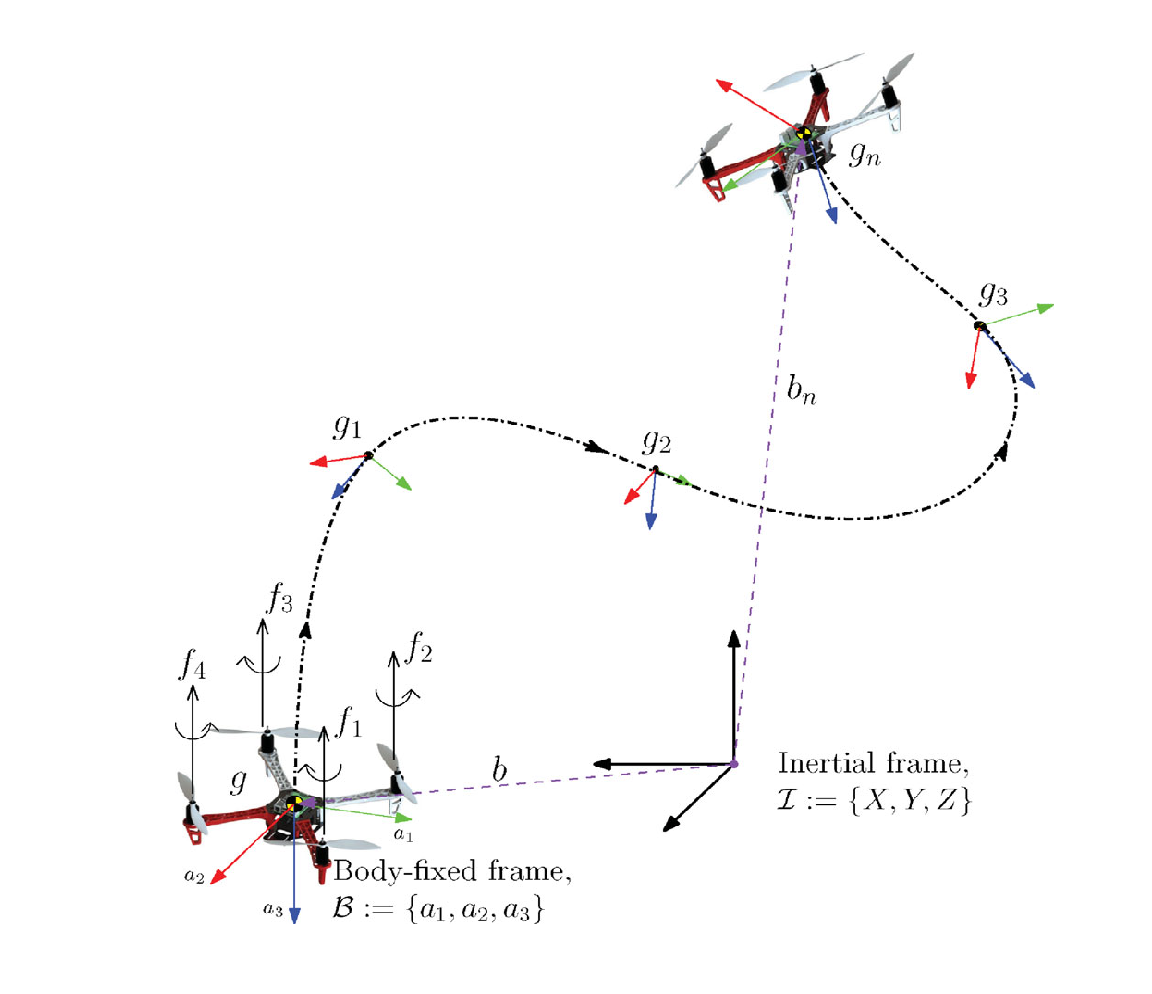}
		\caption{Guidance through a trajectory between initial and final configurations on \SE \cite{hamrah2022finite,viswanathan2018integrated}}
		\label{fig:Configuration}
	\end{figure}

	The configuration of the UAV, modeled as a rigid body, is given by its position and orientation, which are together
	referred to as its pose. To define the pose of the vehicle, we fix a coordinate frame $\mathcal{B}$ to its body and another coordinate frame $\mathcal{I}$ that is fixed in space as the inertial coordinate frame. Define $\textbf{e}_i$ as the unit vector along the axis of the three-dimension space. Let $b \in \bR^3$ denote the position vector of the origin of frame $\mathcal{B}$ with respect to frame $\mathcal{I}$. Let \SO \ denote the orientation (attitude), defined as the rotation matrix from frame  $\mathcal{B}$ to frame  $\mathcal{I}$. The pose of the vehicle can be represented in matrix form as follows:
	\begin{align}\label{eqn:Coordinate 1}
		g = 
		\left[
		\begin{array}{cc}
			R & b \\
			0 & 1 \\
		\end{array}
		\right]
		\in \SE
	\end{align}
	where \SE, the special Euclidean group, is the six-dimensional Lie group of rigid body motions. A diagram of guidance and trajectory tracking on \SE \ through a set of waypoints is presented in Figure \ref{fig:Configuration} as follows.
	\subsection{System kinematics and dynamics}
	The instantaneous pose (position and attitude) is compactly represented by $\mathrm{g} = (b, R) \in  \SE$. The UAV's kinematics is then defined by:
	\begin{align}
		\left\{
		\begin{array}{lr}
			\dot{b} = v = R\nu,\\
			\dot{R} =R\Omega^{\times},
		\end{array}
		\right.
		\label{eqn:Kinematics}
	\end{align}
	where $v \in \R{3}$ and $\nu \in \R{3}$ denote the translational velocity in frames $\mathcal{I}$ and $\mathcal{B}$ respectively, and $\Omega$ $\in$ $ \R{3}$ is the angular velocity in body-fixed frame $\mathcal{B}$. The overall system kinematics and dynamics of a rotorcraft UAV with a body-fixed plane of rotors are given by: 
	\begin{align}\label{eqn:System}
		\begin{cases}
			& \dot{b} = v = R\nu \\
			& m\dot{v} = mg\textbf{e}_3 - fR\textbf{e}_3 + \varphi_D  \\
			& \dot{R} = R\Omega^{\times} \\
			&J\dot{\Omega} = J\Omega\times\Omega + \tau + \tau_D 
		\end{cases}    
	\end{align}
	where $\textbf{e}_3 = [0 \quad 0 \quad 1]\T$, $f \in \bR$ is the scalar thrust force, and $\tau \in \bR^{3}$ is the control torque created by the rotors, $g$ denotes the acceleration due to gravity and $m \in \bR^{+}$ and $J = J\T \in \bR^{3 \times 3}$ are the mass and inertia matrix of the UAV, respectively. The force and torque disturbances are denoted $\varphi_D$ and $\tau_D$ respectively, which are mainly due to unsteady aerodynamics.	
	\subsection{Morse function on \SO}
	The following Lemma is utilized in the rotational ESO and attitude tracking control schemes for the aircraft.
	\begin{Lemma}\label{lem:sK Definition}
		\textup{\cite{bohn2016almost}}
		Consider attitude kinematics 
		\begin{align}\label{eqn:Attitude Kinematics}
			\dot{R}=R\Omega^\times, R\in\SO, \Omega\in \sso.
		\end{align}
		Define $K= \diag(K_1,K_2,K_3)$, where $K_1>K_2>K_3 \ge 1$. Define
		\begin{align}\label{eqn:sK Definition}
			s_K(R)=\sum\limits_{i=1}^{3}K_i(R\Tp \textbf{\textup e}_i)\times \textbf{\textup e}_i, 
		\end{align}
		such that\textup{ $\frac{\di}{\di t}\lan K,I-R\ran=\Omega\T s_K(R)$}. Here \textup{$\lan A, B\ran=\tr  (A\T B)$}, which makes 
		$\lan K,I-R\ran$ a Morse function defined on $\SO$.
		Let $\cS\subset\SO$ be a closed subset containing the identity in its interior, defined by
		\begin{align}\label{eqn:S Definition}
			\begin{split}
				\cS &= \big\{ R\in\SO\, :\, R_{ii}\ge 0 \mbox{\textup{ and} } R_{ij}R_{ji}\le 0, \\
				&\forall i,~ j\in \{1,2,3\},\ i\ne j\big\}. 
			\end{split}
		\end{align}
		Then for $\forall R\in\cS$, we have 
		\begin{align}\label{eqn:sK Bound} 
			s_K(R)\Tp s_K(R) \geq \langle K, I-R \rangle.  
		\end{align}
	\end{Lemma}
	\begin{Remark}[Almost global domain of attraction]\label{rem:Critical Points}
		{\em \cite{sanyal2010almost} We know that the subset of \ \SO \ where $s_K(R) = 0, R \in \SO $, which is also the set of critical points for $\langle I-R, K\rangle$, is 
			\begin{align}\label{eqn:Critical Points}
				\begin{split}
					C &\triangleq \{I, \textup{diag}(1,-1,-1), \textup{diag}(-1,1,-1), \\
					& \textup{diag}(-1,-1,1)\}\subset \SO.
				\end{split}
			\end{align}
			In addition, the global minimum of Morse-Function is $R=I$.
		}
	\end{Remark}
	\subsection{Tracking error kinematics and dynamics}
	Let $g^d(t) \in \SE$ be the desired pose generated by a guidance scheme \cite{viswanathan2018integrated}. Let $v^d$ and $\nu^d$ denote the desired translational velocity in the inertial frame $\mathcal{I}$ and the body frame $\mathcal{B}$, respectively, and $\Omega^d$ denote the body's reference angular velocity in the body frame. Then, the tracking error is given by,
	\begin{align}	\bar h = (\mathrm{g}^d)^{-1}\mathrm{g} =  
		\begin{pmatrix}
			Q & x  \\
			0 & 1 
		\end{pmatrix} \in \SE,
		\label{eqn:Tracking Error}
	\end{align}
	where $Q=(R^d)\T R$ is the attitude tracking error, and $x = (R^d)\T(b-b^d)=(R^d)\T\Tilde{b}$ is the position tracking error, both in the body-fixed frame. Also, the translation velocity tracking error is given by
	\begin{align}\label{eqn:Translational Velocity Error}
		\Tilde{v} = v - v^d, 
	\end{align}
	and the angular velocity tracking error is given by,
	\begin{align}\label{eqn:Attitude Control Error Omega}
		\omega = \Omega - Q\T \Omega^d.
	\end{align}
	Thus, in the inertial frame $\mathcal{I}$, the translational tracking error kinematics and dynamics can be summarized as
	\begin{align}\label{eqn:Translational Tracking Error Kinematics and Dynamics}
		\begin{split}
			\dot{\Tilde{b}} &= \Tilde{v}, \\
			m\dot{\Tilde{v}}&= mg \textbf{e}_3 - fR\textbf{e}_3 + \varphi_D -m\dot{v}^d.
		\end{split}
	\end{align}
	in the body-fixed frame $\mathcal{B}$, the attitude tracking error kinematics and dynamics can be summarized as
	\begin{align}\label{eqn:Attitude Tracking Error Kinematics and Dynamics}
		\begin{split}
			\dot{Q} &= Q\omega^\times, \\
			J\dot{\omega} &= \tau + \tau_D +J(\omega^\times Q\T \Omega^d -Q\T \dot{\Omega}^d )+(J\Omega)\times\Omega . 
		\end{split}
	\end{align}
	The rotational error dynamics is decoupled from the translational error dynamics such that the translation control force, $f$, is obtained in the inertial frame followed by the appropriate attitude tracking control law, $\tau$, in body frame to track the desired trajectory, $b^d$. 
	\subsection{ESO estimations and errors}
	The ESO design on \SE \ is split into translational ESO design on vector space $\bR^3$ and rotational ESO design on \SO. Let ($\widehat{b}, \widehat{v}, \widehat{\varphi}_D )\in \bR^3 \times \bR^3 \times \bR^3$ be the estimated translational position, velocity, and disturbance forces, as the states of translational ESO. The estimation errors of translational ESO are then defined as follows,
	\begin{align}\label{eqn:Translational ESO Error}
		e_b = b - \widehat{b}, e_v = v -\widehat{v}, e_\varphi = \varphi_D - \widehat{\varphi}_D,
	\end{align}
	which are estimation errors of translational position, velocity, and total disturbance force respectively. 
	
	Let $(\widehat{R}, \widehat{\Omega}, \widehat{\tau}_D)\in \SO\times\bR^3\times\bR^3$ be the estimated attitude, angular velocity, and disturbance torque states provided by the rotational ESO. For the rotational ESO, the error states are defined as follows. Please note that the attitude estimation error can be defined as
	\begin{align}\label{eqn:Attitude ESO Error ER}
		E_R = \widehat{R}\T R,
	\end{align}
	on the group of rigid body rotations, $\SO$,  which is not a vector space. 
	The angular velocity estimation error, $e_\Omega$, and torque disturbance estimation error, $e_\tau$, are expressed on the vector space $\bR^3$, and are defined as:
	\begin{align}\label{eqn:Attitude ESO Error}
		e_\Omega = \Omega -E\T_R\widehat{\Omega}, \ e_\tau = \tau_D - \widehat{\tau}_D .
	\end{align}
	
	With a proper ESO design on \SE, the error states $(e_b, e_v, e\varphi)$ and $(E_R, e_\Omega, e_\tau)$ will converge to $(0, 0, 0)$ and $(I, 0, 0)$, respectively. The ESO design and its stability proof will be described in detail in the following section.
	\section{Fast Finite-Time Stable Extended State Observer (FFTS-ESO) on \SE}\label{sec:ESO}
	In this section, we present the FFTS-ESO on \SE. As mentioned in the previous section, the ESO design on \SE \ can be represented as a translational ESO on the vector space $\bR^3$ to estimate disturbance forces, and an rotational ESO on \SO \ to estimate disturbance torques. We present the ESO designs in two results and their stability proofs in this section.  
	\subsection{ESO for Translational Motion}
	\begin{Proposition}[Translational ESO]\label{prop:Translational ESO}
		Consider the following ESO design for the translational motion:
		{\em
			\begin{align}\label{eqn:Translational ESO}
				\begin{split}
					\dot{\widehat{b}} &= \widehat{v}, \\
					m\dot{\widehat{v}} &= mg\textbf{e}_3 -fR\textbf{e}_3 + m k_{t1} \phi_1(\psi_t) \\
					&+ m\kappa_t \Big[ (e\Tp_be_b)^{\frac{1-p}{p}}H\Big(e_b,\frac{p-1}{p}\Big)e_v + e_v\Big] +\widehat{\varphi}_D , \\
					\dot{\widehat{\varphi}}_D &= m k_{t2} \phi_2(\psi_t),
				\end{split}
			\end{align}
		}
		where $\psi_t$ is defined as
		\begin{align}\label{eqn:psit}
			\psi_t = e_v + \kappa_t \Big[e_b + (e\Tp_be_b)^{\frac{1-p}{p}}e_b\Big],
		\end{align}
		and $\phi_1(\cdot)$ is as defined in the expression in \eqref{eqn:phi1phi2}. In addition, we define the constant $k_{t3}$, which does not appear in the expressions \eqref{eqn:Translational ESO} and \eqref{eqn:psit}, but occurs in the terms $\phi_1(\psi_t)$ and $\phi_2(\psi_t)$, where it takes the place of $k_3$ in \eqref{eqn:phi1phi2}.
		The positive scalar gains $k_{t1}, k_{t2}, k_{t3}$, and $\kappa_t$  
		are constrained as follows.
		\begin{itemize}
			\item (Constraint 1) 
			The matrix $\mathcal{A}_t \in \bR^{2\times2}$ defined as:
			\begin{align}\label{eqn:Matrix At}
				\mathcal{A}_t = 
				\begin{bmatrix}
					-k_{t1} & 1 \\ 
					-k_{t2} & 0
				\end{bmatrix},
			\end{align}
			is a Hurwitz matrix.
			\item (Constraint 2) For $\mathcal{A}_t$ as defined above, $\forall\, \mathcal{Q}_t \in \bR^{ 2\times 2}$ where $\mathcal{Q}_t \succ 0$, the Lyapunov equation,
			\begin{align}\label{eqn:Lyapunov Equation TESO 1}
				\mathcal{A}\Tp_t \mathcal{P}_t+\mathcal{P}_t\mathcal{A}_t =-\mathcal{Q}_t,
			\end{align}
			has a unique solution $\mathcal{P}_t$. 
			The eigenvalues of $\mathcal{Q}_t$ are constrained as follows:
			\begin{align}\label{eqn:Translational ESO Constraint}
				k_{t3}\frac{\lambda_\textup{min}\left\{\mathcal{Q}_t\right\}}{ \lambda_\textup{max}\left\{\mathcal{P}_t\right\}} - \frac{1}{k^2_{t3}\lambda_\textup{min}\left\{\mathcal{P}_t\right\}}>0.
			\end{align}
			\item (Constraint 3) $\kappa_t>1/2$.
		\end{itemize}
	\end{Proposition}
	\begin{Theorem}\label{thr:Transaltional ESO Error Dynamics}
		With the observer errors for the translational ESO defined by \eqref{eqn:Translational ESO Error}, the translational kinematics and dynamics given by \eqref{eqn:System}, and the ESO for translational motion given in Proposition \ref{prop:Translational ESO}, the error dynamics of the ESO is given by:
		\begin{align}\label{eqn:Translational ESO Error Dynamics}
			\begin{split}
				\dot{e}_b &= e_v, \\
				m\dot{e}_v &= -mk_{t1} \phi_1(\psi_t) \\
				& - m\kappa_t \Big[ (e\Tp_be_b)^{\frac{1-p}{p}}H\Big(e_b,\frac{p-1}{p}\Big)e_v + e_v\Big] + e_\varphi,\\
				\dot{e}_\varphi &= -mk_{t2}\phi_2(\psi_t)  + \dot{\varphi}_D.
			\end{split}
		\end{align}
		The error dynamics 
		\eqref{eqn:Translational ESO Error Dynamics} is FFTS at the 
		origin $((e_b,e_v,e_\varphi)=(0,0,0))$, when the resultant disturbance force is constant $(\dot{\varphi}_D=0)$, and the gains of the ESO are constrained  according to Proposition \ref{prop:Translational ESO}.
	\end{Theorem}
	\begin{Proof*}
		{\em
			Simplify \eqref{eqn:Translational ESO Error Dynamics} as:
			\begin{align}\label{eqn:Translational ESO Error Dynamics 1}
				\begin{split}
					\dot{\psi}_t &= -k_{t1} \phi_1(\psi_t) + m^{-1} e_\varphi,  \\
					m^{-1}\dot{e}_\varphi &= -k_{t2}\phi_2(\psi_t)  + m^{-1} \dot{\varphi}_D. 
				\end{split}
			\end{align}
			Next, define the Lyapunov function to prove Theorem \ref{thr:Transaltional ESO Error Dynamics}: 
			\begin{align}\label{eqn:Lyapunov Translational ESO}
				V_t &= V_{t0}  + e\Tp_be_b, \mbox{ where } V_{t0} = \zeta\Tp_t \mathcal{P}_t \zeta_t 
			\end{align}
			and $\zeta_t$ is defined as:
			\begin{align*}
				\zeta_t = [\phi\Tp_1(\psi_t), \ m^{-1} e\Tp_\varphi]\Tp.
			\end{align*}
			From Theorem \ref{thr:FFTS Differentiator}, \eqref{eqn:Translational ESO Error Dynamics 1} and \eqref{eqn:Lyapunov Bound 3}, we find that the time-derivative of $V_t$ satisfies:
			\begin{align}\label{eqn:Lyapunov Derivative Translational ESO 1}
				\dot{V}_t &\leq -\gamma_{t1} V_{t0} -\gamma_{t2} V^\frac{1}{p}_{t0} + 2e\T_be_v, 
			\end{align}
			where $\gamma_{t1}$ and $\gamma_{t2}$ are defined by:
			\begin{align}\label{eqn:gammat1 gammat2}
				\begin{split}
					\gamma_{t1} &= k_{t3} \frac{\lambda_\textup{min}\left\{\mathcal{Q}_t \right\}}{\lambda_\textup{max}\left\{\mathcal{P}_t\right\}}, \\
					\gamma_{t2} &= \frac{\lambda_\textup{min}\left\{\mathcal{Q}_t \right\}\lambda_\textup{min}\left\{\mathcal{P}_t\right\}^\frac{p-1}{p}p}{\lambda_\textup{max}\left\{\mathcal{P}_t\right\}(3p-2)}.
				\end{split}
			\end{align}
			Substituting \eqref{eqn:psit} into \eqref{eqn:Lyapunov Derivative Translational ESO 1}, we obtain:
			\begin{align}\label{eqn:Lyapunov Derivative Translational ESO 2}
				\begin{split}
					\dot{V}_t &\leq -\gamma_{t1} V_{t0} -\gamma_{t2} V^\frac{1}{p}_{t0} \\
					&+2 e\Tp_b\Big[\psi_t - \kappa_t e_b- \kappa_t (e\Tp_be_b)^{\frac{1-p}{p}}e_b\Big] \\
					&\leq -\gamma_{t1} V_{t0} -\gamma_{t2} V^\frac{1}{p}_{t0}\\
					& +2 e\Tp_b \psi_t - 2 \kappa_t e\Tp_be_b- 2 \kappa_t (e\T_be_b)^{\frac{1}{p}}  \\
					&\leq -\gamma_{t1} V_{t0} -\gamma_{t2} V^\frac{1}{p}_{t0} \\
					&- 2 \kappa_t e\T_be_b- 2 \kappa_t (e\Tp_be_b)^{\frac{1}{p}} + \psi\Tp_t \psi_t + e\Tp_be_b \\
					&\leq -\Big(\gamma_{t1}-\frac{1}{k^2_{t3}\lambda_\textup{min}\left\{\mathcal{P}_t\right\}}\Big) V_{t0} -\gamma_{t2} V^\frac{1}{p}_{t0} \\
					&- (2\kappa_t-1) e\Tp_be_b- 2 \kappa_t (e\Tp_be_b)^{\frac{1}{p}}. 
				\end{split}
			\end{align}
			Therefore, we further obtain:
			\begin{align}\label{eqn:Lyapunov Derivative Translational ESO 3}
				\dot{V}_t < -\Gamma_{t1} V_t  - \Gamma_{t2} V^{\frac{1}{p}}_t,
			\end{align}
			where 
			\begin{align}\label{eqn:Gammat1 Gammat2}
				\begin{split}
					\Gamma_{t1}  &= \textup{min}\left\{k_{t3} \frac{\lambda_\textup{min}\left\{\mathcal{Q}_t \right\}}{\lambda_\textup{max}\left\{\mathcal{P}_t\right\}} - \frac{1}{k^2_{t3} \lambda_\textup{min}\left\{\mathcal{P}_t\right\}}, 2\kappa_t - 1 \right\}, \\
					\Gamma_{t2}  &=  \textup{min}\left\{\frac{\lambda_\textup{min}\left\{\mathcal{Q}_t \right\}\lambda_\textup{min}\left\{\mathcal{P}_t\right\}^\frac{p-1}{p}p}{\lambda_\textup{max}\left\{\mathcal{P}_t\right\}(3p-2)}, 2\kappa_t \right\} .
				\end{split}
			\end{align}
			Based on \eqref{eqn:Lyapunov Derivative Translational ESO 3}, we conclude that when the resultant disturbance force is constant, and the ESO gains satisfy the constraints 1-3 in Proposition \ref{prop:Translational ESO}, the error dynamics of the ESO  \eqref{eqn:Translational ESO Error Dynamics} is FFTS. This concludes the proof of Theorem \ref{thr:Transaltional ESO Error Dynamics}.
		}
		\qedsymbol{}
	\end{Proof*}
	\subsection{ESO for Rotational Motion}
	\begin{Proposition}[Rotational ESO]\label{prop:Attitude ESO}
		Define $e_R = s_k(E_R)$, where $s_K(\cdot)$ is as defined by Lemma \ref{lem:sK Definition}. 
		Define $e_w(E_R,e_\Omega)$ as follows:
		{\em
			\begin{align}\label{eqn:ew}
				e_w(E_R,e_\Omega) = \frac{\textup{d}}{\textup{d}t} e_R =  \sum_{i=1}^3K_i\textbf{e}_i \times (e_\Omega\times E\T_R \textbf{e}_i).
			\end{align}
		}
		Consider the following ESO design:
		\begin{align}\label{eqn:Attitude ESO}
			\begin{split}
				&\dot{\widehat{R}}= \widehat{R}\widehat{\Omega}^\times, \\
				&\dot{\widehat{\Omega}}= E_RJ^{-1}\Big[J\Omega\times\Omega+\widehat{\tau}_D +\tau + k_{a1} J\phi_1(\psi_a)\Big] \\
				& + E_RJ^{-1}\Big[\kappa_a J (e\Tp_Re_R)^{\frac{1-p}{p}}H\Big(e_R,\frac{p-1}{p}\Big)e_w \Big]\\
				& + E_RJ^{-1}(\kappa_a Je_w)+E_Re^\times_\Omega E\Tp_R\widehat{\Omega}, \\
				&\dot{\widehat{\tau}}_D = J k_{a2} \phi_2(\psi_a), 
			\end{split}
		\end{align}
		where $\psi_a$ is defined as follows: 
		\begin{align}\label{eqn:psia}
			\psi_a = e_\Omega + \kappa_a \Big[e_R + (e_R\Tp e_R)^{\frac{1-p}{p}}e_R\Big].
		\end{align}
		In addition, we define the constant $k_{a3}$, which occurs in the terms $\phi_1(\psi_a)$ and $\phi_2(\psi_a)$,
		where it takes the place of $k_3$ in \eqref{eqn:phi1phi2}. The positive scalar gains $k_{a1}, k_{a2}$, $k_{a3}$, and $\kappa_a$ are constrained as follows.
		\begin{itemize}
			\item (Constraint 1) The matrix $\mathcal{A}_a \in \bR^{2\times2}$ defined as:
			\begin{align}\label{eqn:Matrix Aa}
				\mathcal{A}_a = 
				\begin{bmatrix}
					-k_{a1} & 1 \\ 
					-k_{a2} & 0
				\end{bmatrix},
			\end{align}
			is a Hurwitz matrix.
			\item (Constraint 2) For $\mathcal{A}_a$ as defined above and 
			$\forall\, \mathcal{Q}_a \in \bR^{ 2\times 2}$ where $\mathcal{Q}_a  \succ 0$, the Lyapunov equation: 
			\begin{align}\label{eqn:Lyapunov Equation RESO 1}
				\mathcal{A}\Tp_a \mathcal{P}_a+\mathcal{P}_a\mathcal{A}_a =-\mathcal{Q}_a,
			\end{align}
			has a unique solution $\mathcal{P}_a$. The eigenvalues of $\mathcal{Q}_a$ and $\mathcal{P}_a$ are constrained as follows:
			\begin{align}\label{eqn:Attitude ESO Constraint}
				k_{a3}\frac{\lambda_\textup{min}\left\{\mathcal{Q}_a\right\}}{ \lambda_\textup{max}\left\{\mathcal{P}_a\right\}} - \frac{1}{2k^2_{a3}\lambda_\textup{min}\left\{\mathcal{P}_a\right\}}>0.
			\end{align}
			\item (Constraint 3) $\kappa_a>1/2$.
		\end{itemize}
	\end{Proposition}
	\begin{Theorem}\label{thr:Attitude ESO Error Dynamics}
		With the observer errors for the rotational ESO defined by \eqref{eqn:Attitude ESO Error}, the rotational kinematics and dynamics given by \eqref{eqn:System}, and the ESO for rotational motion given in Proposition \ref{prop:Attitude ESO}, the error dynamics of the ESO is given by:
		\begin{align}\label{eqn:Attitude ESO Error Dynamics}
			\begin{split}
				\dot{E}_R&= E_Re^\times_\Omega, \\
				J\dot{e}_\Omega &= -k_{a1}J \phi_1(\psi_a) \\
				&- \kappa_a J\Big[(e\Tp_Re_R)^{\frac{1-p}{p}}H\Big(e_R,\frac{p-1}{p}\Big)e_w + e_w\Big] + e_\tau,  \\
				\dot{e}_\tau &= -k_{a2}J\phi_2(\psi_a) +\dot{\tau}_D. 
			\end{split}
		\end{align}
		The error dynamics \eqref{eqn:Translational ESO Error Dynamics} is almost globally FFTS (AG-FFTS) at the origin $((E_R,e_\Omega,e_\tau)=(I,0,0))$, when the resultant disturbance torque is constant ($\dot{\tau}_D=0$), and the gains of the ESO are constrained  according to Proposition \ref{prop:Attitude ESO}. 
	\end{Theorem}
	\begin{Proof*}
		{\em
			Simplify \eqref{eqn:Attitude ESO Error Dynamics} as:
			\begin{align}\label{eqn:Attitude ESO Error Dynamics 1}
				\begin{split}
					\dot{\psi}_a &= -k_{a1} \phi_1(\psi_a) + J^{-1} e_\tau, \\
					J^{-1}\dot{e}_\tau &= -k_{a2}\phi_2(\psi_a)+J^{-1}\dot{\tau}_D. 
				\end{split} 
			\end{align}
			Next, define the Morse-Lyapunov function to prove Theorem \ref{thr:Attitude ESO Error Dynamics}: 
			\begin{align}\label{eqn:Lyapunov Attitude ESO}
				\begin{split}
					V_a &= V_{a0}  + \langle K, I-E_R \rangle, \mbox{ where } V_{a0} = \zeta\T_a \mathcal{P}_a \zeta_a
				\end{split}
			\end{align}
			and $\zeta_a$ is defined as:
			\begin{align*}
				\zeta_a = [\phi\T_1(\psi_a), \ J^{-1} e\T_\tau]\T.
			\end{align*}
			From Theorem \ref{thr:FFTS Differentiator}, \eqref{eqn:Attitude ESO Error Dynamics 1} and \eqref{eqn:Lyapunov Bound 3}, { we find that the time-derivative of $V_a$ satisfies:}
			\begin{align}\label{eqn:Lyapunov Derivative Attitude ESO 1}
				\dot{V}_a &\leq -\gamma_{a1} V_{a0} -\gamma_{a2} V^\frac{1}{p}_{a0} + e\T_Re_\Omega,
			\end{align}
			where $\gamma_{a1}$ and $\gamma_{a2}$ are defined by:
			\begin{align}\label{eqn:gammaa1 gammaa2}
				\begin{split}
					\gamma_{a1} &= k_{a3} \frac{\lambda_\textup{min}\left\{\mathcal{Q}_a \right\}}{\lambda_\textup{max}\left\{\mathcal{P}_a\right\}}, \\
					\gamma_{a2} &= \frac{\lambda_\textup{min}\left\{\mathcal{Q}_a \right\}\lambda_\textup{min}\left\{\mathcal{P}_a\right\}^\frac{p-1}{p}p}{\lambda_\textup{max}\left\{\mathcal{P}_a\right\}(3p-2)}.
				\end{split}
			\end{align}
			Substituting \eqref{eqn:psia} into \eqref{eqn:Lyapunov Derivative Attitude ESO 1}, we obtain,
			\begin{align}\label{eqn:Lyapunov Derivative Attitude ESO 2}
				\begin{split}
					\dot{V}_a &\leq -\gamma_{a1} V_{a0} -\gamma_{a2} V^\frac{1}{p}_{a0} \\
					&+ e\Tp_R\Big[\psi_a - \kappa_a e_R- \kappa_a (e\Tp_Re_R)^{\frac{1-p}{p}}e_R\Big] \\
					&\leq -\gamma_{a1} V_{a0} -\gamma_{a2} V^\frac{1}{p}_{a0} \\
					&+ e\Tp_R \psi_a - \kappa_a e\Tp_Re_R- \kappa_a (e\Tp_Re_R)^{\frac{1}{p}}  \\
					&\leq -\gamma_{a1} V_{a0} -\gamma_{a2} V^\frac{1}{p}_{a0} \\
					&- \kappa_a e\Tp_Re_R- \kappa_a (e\Tp_Re_R)^{\frac{1}{p}} + \frac{1}{2}\psi\Tp_a \psi_a + \frac{1}{2} e\Tp_Re_R  \\
					&\leq -\Big(\gamma_{a1}-\frac{1}{2k^2_{a3}\lambda_\textup{min}\left\{\mathcal{P}_a\right\}}\Big) V_{a0} -\gamma_{a2} V^\frac{1}{p}_{a0} \\
					&-\Big(\kappa_a-\frac{1}{2}\Big) e\T_Re_R- \kappa_a (e\T_Re_R)^{\frac{1}{p}} 
				\end{split}
			\end{align}
			By applying Lemma \ref{lem:sK Definition} on \eqref{eqn:Lyapunov Derivative Attitude ESO 1}, we obtain,
			\begin{align}\label{eqn:Lyapunov Derivative Attitude ESO 3}
				\begin{split}
					\dot{V}_a &\leq -\Big(\gamma_{a1}-\frac{1}{2k^2_{a3}\lambda_\textup{min}\left\{\mathcal{P}_a\right\}}\Big) V_{a0} -\gamma_{a2} V^\frac{1}{p}_{a0} \\
					&- \Big(\kappa_a-\frac{1}{2}\Big) \langle K,I-E_R\rangle- \kappa_a \langle K,I-E_R\rangle^{\frac{1}{p}}
				\end{split}
			\end{align}
			Therefore, we further obtain:
			\begin{align}\label{eqn:Lyapunov Derivative Attitude ESO 4}
				\dot{V}_a \leq -\Gamma_{a1} V_a  - \Gamma_{a2} V^{\frac{1}{p}}_a,
			\end{align}
			where, 
			\begin{align}\label{eqn:Gammaa1 Gammaa2}
				\begin{split}
					\Gamma_{a1}&=\textup{min}\left\{k_{a3}\frac{\lambda_\textup{min}\left\{\mathcal{Q}_a \right\}}{\lambda_\textup{max}\left\{\mathcal{P}_a\right\}}-\frac{1}{2k^2_{a3}\lambda_\textup{min}\left\{\mathcal{P}_a\right\}},\kappa_a-\frac{1}{2}\right\},\\
					\Gamma_{a2}&=\textup{min}\left\{\frac{\lambda_\textup{min}\left\{\mathcal{Q}_a \right\}\lambda_\textup{min}\left\{\mathcal{P}_a\right\}^\frac{p-1}{p}p}{\lambda_\textup{max}\left\{\mathcal{P}_a\right\}(3p-2)}, \kappa_a \right\}.
				\end{split}
			\end{align}
			Consider the expression given by \eqref{eqn:Lyapunov Derivative Attitude ESO 4}, the set where $\dot{V}_a=0$ is:
			\begin{align}\label{eqn:Attitude ESO Equilibrium Point Set 1}
				\begin{split}
					&\dot{V}_a^{-1}(0)=\left\{(E_R,e_\Omega,e_\tau): s_K(E_R)=0, \mbox{and} \ \zeta_a = 0 \right\} \\
					&= \left\{(E_R,e_\Omega,e_\tau): E_R \in C, e_\Omega=0, \mbox{and} \ e_\tau=0 \right\},
				\end{split}
			\end{align}
			where $C$ is defined by \eqref{eqn:Critical Points}, which express the set of the critical points. With Theorem 8.4 from Khalil \cite{khalil2002nonlinear}, we conclude that $(E_R,e_\Omega,e_\tau)$ converge to the set:
			\begin{align}\label{eqn:Attitude ESO Equilibrium Point Set 2}
				\begin{split}
					&S=\left\{(E_R,e_\Omega,e_\tau)\in \SO\times\bR^3\times\bR^3:  \right.\\
					&\left.E_R \in C, e_\Omega=0, \mbox{and} \ e_\tau=0 \right\},
				\end{split}
			\end{align}
			in finite time.
			Based on \eqref{eqn:Lyapunov Derivative Attitude ESO 4}, and Lemma \ref{lem:FFTS}, we conclude that when the ESO gains satisfy the constraints in Proposition \ref{prop:Attitude ESO},  the set of equilibrium $S$ for the error dynamics \eqref{eqn:Attitude ESO Error Dynamics} is fast finite time stable. 
			
			In $S$, the only stable equilibrium is $(I,0,0)$, while the other three are unstable. The resulting closed-loop system with the estimation errors gives rise to a H\"{o}lder-continuous feedback with exponent less than one $(1/2<1/p<1)$, while in the limiting case of $p=1$, the feedback system is Lipschitz-continuous. Proceeding with a topological equivalence-based analysis similar to the one by Bohn et al.\cite{bohn2016almost}, we conclude that the equilibrium and the
			corresponding regions of attraction of the rotational ESO with $p\in]1,2[$ are identical to those of the corresponding Lipschitz-continuous asymptotically stable ESO with $p=1$, and the region of attraction is almost global. 
			
			To summarize, we conclude that the error dynamics \eqref{eqn:Translational ESO Error Dynamics} is almost globally FFTS (AG-FFTS) at the origin $((E_R,e_\Omega,e_\tau)=(I,0,0))$ when the resultant disturbance torque is constant ($\dot{\tau}_D=0$), and the gains of the ESO are constrained  according to Proposition \ref{prop:Attitude ESO}. This concludes the proof of Theorem \ref{thr:Attitude ESO Error Dynamics}. 
		}
		\qedsymbol{}
	\end{Proof*}
	\begin{Remark}[Disturbance robustness of the ESO]
		{\em
			Consider Corollary \ref{cor:FTS Differentiator Disturbance Robustness} and its constraints on differentiator gains. When the disturbance forces and torques are time-varying, then $\|\dot{\varphi}_D\|, \|\dot{\tau}_D\|>0$. Further, if the constraints on gains in Corollary \ref{cor:FTS Differentiator Disturbance Robustness} are fulfilled, the estimation error dynamics of the proposed ESO will be practically finite-time stable (PFTS). 
		}
	\end{Remark}
	\begin{Remark}[Noise robustness of the ESO]\label{rem:Noise}
		{\em
			Consider Corollary \ref{cor:FTS Differentiator Noise Robustness} and its constraints on differentiator gains. When the ESO measurements have noise and the constraints on gains in Corollary \ref{cor:FTS Differentiator Noise Robustness} are fulfilled, the estimation error dynamics of the proposed ESO will be PFTS. Moreover, according to Lemma \ref{lem:PFTS} and Corollary \ref{cor:FTS Differentiator Noise Robustness}, the $\eta$ in \eqref{eqn:Preliminary PFTS Lyapunov Inequality} of Lemma \ref{lem:PFTS} is a function on the level of noise in information on $R$, $\Omega$, $b$ and $v$ and is monotonically increasing with the level of noise.
		}
	\end{Remark}
	\begin{Remark}[Comparative Analysis of Noise Robustness: FFTS-ESO vs. the FxTSDO by Liu et al.~\cite{liu2022fixed} ]\label{rem:Compare}
		{\em
			We investigate the disturbance (forces or torques) observers proposed by Liu et al.~\cite{liu2022fixed} in their Theorems 1 and 2, known as FxTSDO. The input of FxTSDO relies on the motion signals, $X_2$, $Y_2$,  which represent translational and angular velocities, and $\dot{X}_2$, $\dot{Y}_2$, which represent translational and angular accelerations, respectively. 
			
			However, the high-level noise associated with the translational acceleration obtained from an accelerometer restricts its direct use in a flight control scheme. Additionally, direct measurement of angular acceleration is usually not feasible. 
			
			Furthermore, if $\dot{X}_2$ and $\dot{Y}_2$ are obtained from the finite difference of $X_2$ and $Y_2$, they will have higher noise levels than $X_2$ and $Y_2$, leading to inferior disturbance estimation performance. 
			
			In contrast to FxTSDO, the proposed FFTS-ESO incorporates position and attitude signals, which are zero-order derivatives of motions with lower noise levels. Consequently, FFTS-ESO outperforms FxTSDO in terms of disturbance estimation performance, despite the theoretical fixed-time stability of FxTSDO. We show this through our numerical simulations in Section \ref{sec:Numerical}.
		}
	\end{Remark}
	\section{Fast-Finite Time Stable Active Disturbance Rejection Control (FFTS-ADRC) on \SE}\label{sec:Tracking}
	A robust ADRC on \SE \ is presented in this section and it is split into position and attitude tracking modules. The proposed FFTS-ESO on \SE \ presented in Section \ref{sec:ESO} is utilized here to provide disturbance estimates $\widehat{\varphi}_D$ and $\widehat{\tau}_D$. For tracking control, $e_\varphi$, and $e_\tau$ are not only disturbance estimation errors in \eqref{eqn:Translational ESO Error Dynamics} and \eqref{eqn:Attitude ESO Error Dynamics}, {but they are also the disturbance rejection errors}. The stability proof of the proposed ADRC includes both tracking error dynamics in \eqref{eqn:Translational Tracking Error Kinematics and Dynamics}-\eqref{eqn:Attitude Tracking Error Kinematics and Dynamics} and ESO estimation error dynamics in \eqref{eqn:Translational ESO Error Dynamics} and \eqref{eqn:Attitude ESO Error Dynamics}. This stability proof does not treat the disturbance estimation/rejection errors $e_\varphi$ and $e_\tau$ as zero vectors but as error terms updated according to propositions \ref{prop:Translational ESO} and \eqref{prop:Attitude ESO}. This means that the 
	disturbance estimation errors in the proposed ESO are designed to converge much faster than the tracking errors in ADRC. With our Lyapunov stability analysis, we analyze the coupling between the tracking control and ESO schemes, and how that influences gain tuning for the ADRC scheme. 
	\subsection{ADRC for Translational Motion Control}
	\begin{Proposition}[Position Tracking ADRC]\label{prop:Translational Tracking}
		Given the tracking error kinematics and dynamics in \eqref{eqn:Translational Tracking Error Kinematics and Dynamics}, consider the translational motion tracking control law:	
		{\em
			\begin{align}\label{eqn:Translation Tracking}
				\begin{split}
					\varphi &=fR\textbf{e}_3 = mg\textbf{e}_3 + k_{TD}L_T\Big[\psi_T+(\psi\Tp_T\psi_T)^{\frac{1-p}{p}}\psi_T \Big]\\
					&+ k_{TP}L_T\Tilde{b}  + m\kappa_{T}\Big[\Tilde{v}+(\Tilde{b}\Tp\Tilde{b})^{\frac{1-p}{p}}H\Big(\Tilde{b},\frac{p-1}{p}\Big)\Tilde{v}\Big]\\
					& -m\dot{v}_d +\widehat{\varphi}_D,
				\end{split}
			\end{align}
		}
		where $\psi_T$ is defined as:
		\begin{align}\label{eqn:psiT}
			\psi_T = \Tilde{v} + \kappa_{T} \Big[\Tilde{b} + (\Tilde{b}\Tp\Tilde{b})^{\frac{1-p}{p}}\Tilde{b}\Big],
		\end{align} 
		and $\widehat{\varphi}_D$ is obtained from the translational ESO in Proposition \ref{prop:Translational ESO}. In addition to the ESO gains defined in Proposition \ref{prop:Translational ESO}, define positive scalar control gains $\kappa_{T}$, $k_{TD}$, $k_{TP}$, and a positive definite diagonal matrix $L_T = \textup{diag}(L_{T1},L_{T2},L_{T3})$ that satisfy the following constraints:
		
		(Constraint 1) $k_{TD}$ and $L_T$ are constrained as: 
		\begin{align*}
			k_{TD}\lambda_{\textup{min}}\{L_T\}-\frac{1}{2}>0.
		\end{align*}
		
		(Constraint 2) The decay constant for the translational ESO $\Gamma_{t1}$ defined by \eqref{eqn:Gammat1 Gammat2}, $\mathcal{P}_t$ given by \eqref{eqn:Lyapunov Equation TESO 1}, and the rotorcraft mass $m$ are constrained as:
		\begin{align*} 
			\Gamma_{t1}-\frac{m^2}{2\lambda_{\textup{min}}\{\mathcal{P}_t\}}> 0.
		\end{align*}
	\end{Proposition}
	\begin{Theorem}\label{thr:Translational Tracking Error Dynamics}
		With the translational tracking errors defined by \eqref{eqn:Tracking Error}, kinematics and dynamics given by \eqref{eqn:Translational Tracking Error Kinematics and Dynamics}, ESO error dynamics given by \eqref{eqn:Translational ESO Error Dynamics}, and the ADRC for position tracking given in Proposition \ref{prop:Translational Tracking}, the translational tracking error dynamics satisfies:
		\begin{align}\label{eqn:Translational Tracking Error Dynamics}
			\begin{split}
				\dot{\Tilde{b}} & = \Tilde{v} \\
				m\dot{\Tilde{v}} &= e_\varphi -  k_{TD}L_T\Big[\psi_T+(\psi\Tp_T\psi_T)^{\frac{1-p}{p}}\psi_T \Big]\\
				& - k_{TP}L_T\Tilde{b}  - m\kappa_{T}\Big[\Tilde{v}+(\Tilde{b}\Tp\Tilde{b})^{\frac{1-p}{p}}H\Big(\Tilde{b},\frac{p-1}{p}\Big)\Tilde{v}\Big],
			\end{split}
		\end{align}
		where $e_\varphi$ is updated by the translational ESO error dynamics according to \eqref{eqn:Translational ESO Error Dynamics}. The tracking error dynamics \eqref{eqn:Translational Tracking Error Dynamics} combined with the ESO error dynamics \eqref{eqn:Translational ESO Error Dynamics} is FFTS at the origin $((\Tilde{b},\Tilde{v},e_b,e_v,e_\varphi)=(0,0,0,0,0))$, when the resultant disturbance force is constant $(\dot{\varphi}_D=0)$, and the gains of the ESO and ADRC schemes are constrained according to Proposition \ref{prop:Translational ESO} and Proposition \ref{prop:Translational Tracking}.
	\end{Theorem}
	\begin{Proof*}
		{\em
			Simplify \eqref{eqn:Translational Tracking Error Dynamics} as:
			\begin{align}\label{eqn:Translational Tracking Error Dynamics 1}
				\begin{split}
					\dot{\Tilde{b}} & = \Tilde{v} \\
					m\dot{\psi}_T &= e_\varphi -  k_{TD}L_T\Big[\psi_T+(\psi\T_T\psi_T)^{\frac{1-p}{p}}\psi_T \Big] \\
					&- k_{TP}L_T\Tilde{b} .
				\end{split}
			\end{align}
			Next, define the Lyapunov function:		
			\begin{align}\label{eqn:Translational ADRC Lyapunov Function}
				V_T =V_t+\frac{1}{2}m\psi\T_T \psi_T + \frac{1}{2}k_{TP}\Tilde{b}\T\Tilde{b}, 
			\end{align}			
			where $V_t$ is defined by \eqref{eqn:Lyapunov Translational ESO} in the proof of Theorem \ref{thr:Transaltional ESO Error Dynamics}. We obtain the time-derivative of this Lyapunov function \eqref{eqn:Translational ADRC Lyapunov Function} as:			
			\begin{align}\label{eqn:Translational ADRC Lyapunov Derivative 1}
				\begin{split}
					\dot{V}_T &= m\psi\Tp_T \dot{\psi}_T+ k_{TP} \Tilde{b}\Tp\dot{\Tilde{b}} + \dot{V}_t\\
					&\leq \psi\Tp_T e_\varphi - k_{TD} \psi\T_T L_{T} \psi_T\\
					&- k_{TD} (\psi\Tp_T\psi_T)^{\frac{1-p}{p}} \psi\Tp_T L_{T} \psi_T\\
					&-k_{TP}\psi_T\Tp \Tilde{b} +k_{TP}\Tilde{b}\Tp \Tilde{v} -\Gamma_{t1}V_t - \Gamma_{t2}V_t^\frac{1}{p} \\
					&\leq \frac{1}{2}\psi\T_T\psi_T +\frac{1}{2} e_\varphi\T e_\varphi \\
					&-k_{TD}\lambda_{\textup{min}}\{L_T\}\psi\T_T \psi_T -k_{TD}\lambda_{\textup{min}}\{L_T\}(\psi\T_T \psi_T)^\frac{1}{p} \\
					&-k_{TP}\psi_T\T \Tilde{b} +k_{TP}\Tilde{b}\T \Big[\psi_T - \kappa_T(\Tilde{b}+(\Tilde{b}\T \Tilde{b})^{\frac{1-p}{p}}\Tilde{b}) \Big] \\
					&-\Gamma_{t1}V_t - \Gamma_{t2}V_t^\frac{1}{p}  
				\end{split}
			\end{align}
			\begin{align*}
				\begin{split}
					&\leq -\Big(k_{TD}\lambda_{\textup{min}}\{L_T\}-\frac{1}{2}\Big)\psi\Tp_T \psi_T \\
					&-k_{TD}\lambda_{\textup{min}}\{L_T\}(\psi\Tp_T \psi_T)^\frac{1}{p}\\
					&-\kappa_T k_{TP}(\Tilde{b}\T \Tilde{b})- \kappa_T k_{TP}(\Tilde{b}\T \Tilde{b})^{\frac{1}{p}}\\
					&-\Big(\Gamma_{t1}-\frac{m^2}{2\lambda_{\textup{min}}\{\mathcal{P}_t\}}\Big)V_t - \Gamma_{t2}V_t^\frac{1}{p} \\
					&\leq -\left(2k_{TD}\lambda_{\textup{min}}\{L_T\}-1\right)m^{-1}\Big(\frac{1}{2}m\psi\Tp_T \psi_T\Big) \\
					&-2^\frac{1}{p}k_{TD}\lambda_{\textup{min}}\{L_T\}m^{-\frac{1}{p}}\left(\frac{1}{2}m\psi\Tp_T \psi_T\right)^\frac{1}{p} \\
					&-2\kappa_T\Big(\frac{1}{2} k_{TP}\Tilde{b}\Tp \Tilde{b} \Big)- 2^\frac{1}{p}\kappa_T k_{TP}^\frac{p-1}{p}\Big(\frac{1}{2} k_{TP}\Tilde{b}\T \Tilde{b} \Big)^{\frac{1}{p}} \\
					&-\Big(\Gamma_{t1}-\frac{m^2}{2\lambda_{\textup{min}}\{\mathcal{P}_t\}}\Big)V_t - \Gamma_{t2}V_t^\frac{1}{p},  
				\end{split}
			\end{align*}
			where $\Gamma_{t1}$, $\Gamma_{t2}$ and $\mathcal{P}_t$ are as defined in the proof of Theorem \ref{thr:Transaltional ESO Error Dynamics}. Thus, we obtain the following inequality:
			\begin{align}\label{eqn:Translational ADRC Lyapunov Derivative 2}
				\dot{V}_T \leq -\Gamma_{T1}V_T -\Gamma_{T2}V^\frac{1}{p}_T,
			\end{align}
			where $\Gamma_{T1}$ and $\Gamma_{T2}$ are defined by:
			\begin{align*}
				\begin{split}
					&\Gamma_{T1}  =\\
					&\textup{min}\left\{\left(2k_{TD}\lambda_{\textup{min}}\{L_T\}-1\right)m^{-1}, 2\kappa_{T}, \Gamma_{t1}-\frac{m^2}{2\lambda_{\textup{min}}\{\mathcal{P}_t\}} \right\}, \\ 
					&\Gamma_{T2}  =\\
					&\textup{min}\left\{2^\frac{1}{p}k_{TD}\lambda_{\textup{min}}\{L_T\}m^{-\frac{1}{p}},  2^\frac{1}{p}\kappa_T k_{TP}^\frac{p-1}{p}, \Gamma_{t2} \right\}. 
				\end{split}
			\end{align*}
			Based on \eqref{eqn:Translational ADRC Lyapunov Derivative 2}, we conclude that when the resultant disturbance force is constant ($\dot{\varphi}_D = 0$), and the ESO and ADRC gains satisfy the constraints in Proposition \ref{prop:Translational ESO} and Proposition \ref{prop:Translational Tracking}, the tracking error dynamics  \eqref{eqn:Translational Tracking Error Dynamics} for the translational motion is FFTS. This concludes the proof of Theorem \ref{thr:Translational Tracking Error Dynamics}.
		}\qedsymbol{}   
	\end{Proof*}
	With this translational motion control scheme, a desired control force vector $\varphi=fR\textbf{e}_3\in \bR^3$ is generated. With $fR\textbf{e}_3$ obtained in this manner, the methodology utilized in \S 3.3 of \cite{viswanathan2018integrated} can be employed to generate the desired (reference) attitude profile to be tracked by the attitude control system, which is described in the following subsection. The term $r_{3d}$ in \cite{viswanathan2018integrated}, which denotes the third column of the rotation matrix, is re-defined here as $ r_{3d}:= {\varphi}/{||\varphi||}.$ The rest of the tracking control design is similar to what has already been used in our prior research \cite{viswanathan2018integrated,hamrah2022finite}. 
	
	\subsection{ADRC for Rotational Motion Control}
	\begin{Proposition}[Attitude Tracking ADRC]\label{prop:Attitude Tracking}
		Given the tracking error kinematics and dynamics in \eqref{eqn:Attitude Tracking Error Kinematics and Dynamics}, consider the attitude tracking control law
		{\em
			\begin{align}\label{eqn:Attitude Tracking}
				\begin{split}
					&\tau=- {k_{AD}L_A}\Big[\psi_A+(\psi\Tp_A\psi_A)^{\frac{1-p}{p}}\psi_A\Big] -k_{AP}s_K(Q)\\
					&-k_{AI}\psi_{AI}   -J(Q\Tp {\dot{\Omega}}^d-\omega^{\times} Q\Tp \Omega^d)-J\Omega\times \Omega- \widehat{\tau}_D\\
					&-\kappa_{A}Jw(Q,\omega)\\
					&-\kappa_{A}J(s_K(Q)\Tp s_K(Q))^{\frac{1-p}{p}}H\Big(s_K(Q),\frac{p-1}{p}\Big)w(Q,\omega), \\
					&{ \dot{\psi}_{AI}} {= -L_A \psi_{AI} -  L_A(\psi_{AI}\Tp\psi_{AI})^{\frac{1-p}{p}} \psi_{AI} + \psi_A  }
				\end{split}
			\end{align}
		}
		where $\psi_{AI}$ is defined as an integral term initialized with $\psi_{AI}(0)=0$, $\psi_A$ is defined as: 	
		{\em	
			\begin{equation}\label{eqn:psiA}
				\psi_A = \omega + \kappa_{A} \Big[s_K(Q) + (s_K(Q)\Tp s_K(Q) )^{\frac{1-p}{p}}s_K(Q) \Big],
			\end{equation}
		} 
		$s_K(Q)$ is defined by Lemma \ref{lem:sK Definition}, $w(Q,\omega)$ is defined as
		{\em
			\begin{align}\label{eqn:w}
				w(Q,\omega) = \frac{\textup{d}}{\textup{d}t} s_K(Q) =  \sum_{i=1}^3K_i\textbf{e}_i \times (\omega \times Q\Tp \textbf{e}_i),
			\end{align}
		}
		and $\widehat{\tau}_D$ is obtained from the rotational ESO in Proposition \ref{prop:Attitude ESO}.  In addition to the ESO gains defined in Proposition \ref{prop:Attitude ESO}, define positive scalar gains $\kappa_{A}$, $k_{AD}$, $k_{AP}$, $k_{AI}$, and a positive definite diagonal matrix gain $L_A = \textup{diag}(L_{A1},L_{A2},L_{A3})$ that satisfy the following constraints:
		
		(Constraint 1) $k_{AD}$ and $L_A$ are constrained as $$2 k_{AD} \lambda_{\textup{min}}\{L_A\}- 1 > 0. $$
		
		(Constraint 2) The decay constant for the rotational ESO $\Gamma_{a1}$ defined by \eqref{eqn:Gammaa1 Gammaa2}, $\mathcal{P}_a$ given by \eqref{eqn:Lyapunov Equation RESO 1}, and the rotorcraft inertia $J$ are constrained as:    
		$$\Gamma_{a1} - \frac{1}{2} \lambda_{\textup{min}}\{J^{-2}\}^{-1}\lambda_{\textup{min}}\{\mathcal{P}_a\}^{-1}> 0.$$
		
	\end{Proposition}
	\begin{Theorem}\label{thr:Attitude Tracking Error Dynamics}
		With the attitude tracking errors defined by \eqref{eqn:Tracking Error}, kinematics and dynamics given by \eqref{eqn:Attitude Tracking Error Kinematics and Dynamics}, ESO error dynamics given by \eqref{eqn:Attitude ESO Error Dynamics}, and the ADRC for attitude tracking given in Proposition \ref{prop:Attitude Tracking}, the rotational tracking error dynamics of the ADRC is given by:
		{\em
			\begin{align}\label{eqn:Attitude Tracking Error Dynamics}
				\begin{split}
					\dot{Q} & = Q\omega \\
					J\dot{\omega} &= e_\tau -{k_{AD}L_A}\Big[\psi_A+(\psi\T_A\psi_A)^{\frac{1-p}{p}}\psi_A\Big] \\
					&-k_{AP}s_K(Q)-k_{AI}\psi_{AI}-\kappa_AJw(Q,\omega)  \\ 
					&-\kappa_AJ(s_K(Q)\T s_K(Q) )^{\frac{1-p}{p}}H\Big(s_K(Q),\frac{p-1}{p}\Big)w(Q,\omega) \\
					{\dot{\psi}_{AI}} &{ = -L_A \psi_{AI} -  L_A(\psi_{AI}\Tp\psi_{AI})^{\frac{1-p}{p}} \psi_{AI} + \psi_A,} 
				\end{split}
			\end{align}
		}
		where $e_\tau$ is updated by the rotational ESO error dynamics given by  \eqref{eqn:Attitude ESO Error Dynamics}. The tracking error dynamics given by \eqref{eqn:Attitude Tracking Error Dynamics}, combined with the ESO error dynamics given by \eqref{eqn:Attitude ESO Error Dynamics}, is almost globally fast finite-time stable (AG-FFTS) at  $((Q,\omega,\psi_{AI}, E_R, e_\Omega, e_\tau)=(I,0,0,I,0,0))$, when the resultant disturbance torque is constant ($\dot{\tau}_D=0$), and the gains of the ESO and ADRC schemes are constrained according to Proposition \ref{prop:Attitude ESO} and Proposition \ref{prop:Attitude Tracking}.
	\end{Theorem}
	\begin{Proof*}
		{\em
			Simplify \eqref{eqn:Attitude Tracking Error Dynamics} as:
			\begin{align}\label{eqn:dot psiA}
				\begin{split}
					\dot{Q} & = Q\omega \\
					J\dot{\psi}_A &= e_\tau -k_{AD}L_A\Big[\psi_A+(\psi\Tp_A\psi_A)^{\frac{1-p}{p}}\psi_A\Big] \\
					&-k_{AP}s_K(Q) - k_{AI}\psi_{AI} \\
					\dot{\psi}_{AI} & = -L_A \psi_{AI} -  L_A(\psi_{AI}\Tp\psi_{AI})^{\frac{1-p}{p}} \psi_{AI} + \psi_A   . 
				\end{split}
			\end{align}
			Next, define the following Morse-Lyapunov function:		
			\begin{align}\label{eqn:Attitude ADRC Lyapunov Function}
				\begin{split}
					V_A &= V_a+\frac{1}{2}\psi\T_A J\psi_A + k_{AP}\langle K, I-Q\rangle \\
					&+ \frac{1}{2} k_{AI} \psi\T_{AI} \psi_{AI},
				\end{split}
			\end{align}
			where $V_a$ is defined by \eqref{eqn:Lyapunov Attitude ESO} in the proof of Theorem \ref{thr:Attitude ESO Error Dynamics}. We obtain the time derivative of $V_A$ as:
			\begin{align}\label{eqn:Attitude ADRC Lyapunov Derivative 1}
				\begin{split}
					\dot{V}_A &= \dot{V}_a+ \psi\T_A J\dot{\psi}_A + k_{AP}s_K(Q)\T \omega + k_{AI}\psi\T_{AI} \dot{\psi}_{AI} \\
					&\leq -\Gamma_{a1} V_a  - \Gamma_{a2} V^{\frac{1}{p}}_a+ \psi\T_A e_\tau \\
					&-k_{AD}\psi\T_AL_{A}\Big[\psi_A+(\psi\T_A\psi_A)^{\frac{1-p}{p}}\psi_A\Big] + k_{AI}\psi_{AI}\T \psi_A\\
					&-k_{AP} \psi\T_A s_K(Q) -k_{AI} \psi\T_A \psi_{AI} + k_{AP} \psi\T_A s_K(Q)    \\
					&-k_{AP}\kappa_{A}s_K(Q)\Tp s_K(Q) -k_{AP}\kappa_{A}(s_K(Q)\T s_K(Q) )^{\frac{1}{p}}  \\
					&-k_{AI}\psi_{AI}\T L_A \psi_{AI} -  k_{AI}(\psi_{AI}\T\psi_{AI})^{\frac{1-p}{p}} \psi_{AI}\T L_A \psi_{AI}  \\
					&\leq -\Gamma_{a1} V_a  - \Gamma_{a2} V^{\frac{1}{p}}_a + \frac{1}{2}e_\tau\T e_\tau \\
					&-\Big(k_{AD} \lambda_{\textup{min}}\{L_A\}-\frac{1}{2}\Big) \psi\T_A \psi_A \\
					&-k_{AD} \lambda_{\textup{min}}\{L_A\}(\psi\T_A \psi_A)^\frac{1}{p} \\
					&-k_{AP}\kappa_A \langle K, I-Q\rangle - k_{AP}\kappa_A\langle K, I-Q\rangle^\frac{1}{p} \\
					&-k_{AI}\lambda_{\textup{min}}\{L_A\}\psi_{AI}\T \psi_{AI} -  k_{AI}\lambda_{\textup{min}}\{L_A\}(\psi_{AI}\T\psi_{AI})^{\frac{1}{p}} \\
					&\leq -\Big[\Gamma_{a1} - \left(2 \lambda_{\textup{min}}\{J^{-2}\}\lambda_{\textup{min}}\{\mathcal{P}_a\}\right)^{-1}\Big] V_a - \Gamma_{a2} V^{\frac{1}{p}}_a \\
					&-\kappa_A k_{AP}\langle K, I-Q\rangle - \kappa_Ak_{AP}^{\frac{p-1}{p}} (k_{AP} \langle K, I-Q\rangle)^\frac{1}{p} \\
					&-\Big(2k_{AD} \lambda_{\textup{min}}\{L_A\}-1\Big)\lambda_{\textup{max}}\{J\}^{-1} \Big(\frac{1}{2}\psi\T_A J \psi_A \Big)\\
					&-2^{\frac{1}{p}} k_{AD} \lambda_{\textup{min}}\{L_A\}\lambda_{\textup{max}}\{J\}^{-\frac{1}{p}}\Big(\frac{1}{2}\psi\T_A J \psi_A \Big)^\frac{1}{p}  \\
					&-2\lambda_{\textup{min}}\{L_A\} \Big(\frac{1}{2}k_{AI}\psi_{AI}\T \psi_{AI}\Big) \\
					&-2^{\frac{1}{p}} k_{AP}^{\frac{p-1}{p}}\lambda_{\textup{min}}\{L_A\}\Big(\frac{1}{2}k_{AI}\psi_{AI}\T \psi_{AI}\Big)^{\frac{1}{p}},
				\end{split}
			\end{align}
			where $\Gamma_{a1}$, $\Gamma_{a2}$ and $\mathcal{P}_a$ are defined by \ref{eqn:gammaa1 gammaa2} in the proof of Theorem \ref{thr:Attitude ESO Error Dynamics}. Simplify \eqref{eqn:Attitude ADRC Lyapunov Derivative 1}:
			\begin{align}\label{eqn:Attitude ADRC Lyapunov Derivative 2}
				\dot{V}_A   \leq -\Gamma_{A1}V_A -\Gamma_{A2}V^\frac{1}{p}_A
			\end{align}
			where $\Gamma_{A1}$ and $\Gamma_{A2}$ are defined as:
			\begin{align*}
				\Gamma_{A1}  &= \textup{min}\Bigg\{
				\Big(\Gamma_{a1} - \frac{1}{2} \lambda_{\textup{min}}\{J^{-2}\}^{-1}\lambda_{\textup{min}}\{\mathcal{P}_a\}^{-1}\Big), \Bigg. \\
				&\Bigg. 2\Big(k_{AD} \lambda_{\textup{min}}\{L_A\}-\frac{1}{2}\Big)\lambda_{\textup{max}}\{J\}^{-1}, \
				\kappa_A, \
				2 \lambda_{\textup{min}}\{L_A\} 
				\Bigg\},  \\
				\Gamma_{A2}  &= \textup{min}\Bigg\{
				\Gamma_{a2},  \
				2^{\frac{1}{p}} k_{AD} \lambda_{\textup{min}}\{L_A\}\lambda_{\textup{max}}\{J\}^{-\frac{1}{p}}, \Bigg. \\
				&\Bigg. \kappa_Ak_{AP}^{\frac{p-1}{p}}, \
				2^{\frac{1}{p}} k_{AP}^{\frac{p-1}{p}}\lambda_{\textup{min}}\{L_A\}
				\Bigg\}. \notag
			\end{align*}
			Based on \eqref{eqn:Attitude ADRC Lyapunov Derivative 2}, we conclude that when the resultant disturbance torque is constant ($\dot{\tau}_D = 0$), and the ESO and ADRC gains satisfy the constraints in Proposition \ref{prop:Attitude ESO} and Proposition \ref{prop:Attitude Tracking}, the tracking error dynamics \eqref{eqn:Attitude Tracking Error Dynamics} is AG-FFTS. Since the discussion on the equilibrium set including \SO \ is covered in the proof of Theorem \ref{thr:Attitude ESO Error Dynamics}, we omit it here in this proof for brevity. This concludes the proof of Theorem \ref{thr:Attitude Tracking Error Dynamics}.
		}\qedsymbol{}   
	\end{Proof*}
	\section{Numerical Simulations}\label{sec:Numerical}
	In this section, we present several sets of numerical simulations, arranged in two subsections to validate the proposed ESO and ADRC. In subsection \ref{subsec:comparison}, we compare the proposed FFTS-ESO with the existing ESO \cite{shao2018robust} and DO \cite{liu2022fixed} designs, on their disturbance estimation performance in different flight scenarios. Subsection \ref{subsec:LGVI} describes the numerical simulation of the ADRC proposed in Section \ref{sec:Tracking} in different flight scenarios with the controlled dynamics numerically simulated by a geometric integrator\cite{nordkvist2010lie}.
	\vspace{-3mm}
	\subsection{The comparison between FFTS-ESO and other ESO/DO}\label{subsec:comparison}
	We compare the proposed FFTS-ESO with existing disturbance estimation schemes, which are LESO \cite{shao2018robust} and FxTSDO \cite{liu2022fixed}, on their disturbance estimation performance in four different simulated flight scenarios, with and without the presence of measurement noises. The four flight scenarios correspond to four desired trajectories. The inertia and mass of the simulated rotorcraft UAV are $  J=\textup{diag}([0.0820,0.0845,0.1377]) \ \textup{kg}\cdot \textup{m}^2, \quad m = 4.34 \ \textup{kg}$\cite{pounds2010modelling}. Since the target of the simulation is to validate and compare the disturbance estimation performance, the actuator dynamics and saturation are not included in the results reported in this section. The tracking control scheme to drive the UAV to track the desired trajectories is the tracking control scheme reported in Section \ref{sec:Tracking} without disturbance rejection terms, such that $\widehat{\tau}_D=0$ in \eqref{eqn:Attitude Tracking} and $\widehat{\varphi}_D=0$ in \eqref{eqn:Translation Tracking}. 
	We use the MATLAB/Simulink software with ODE4 (Runge-Kutta fourth order) solver to conduct this set of simulations. The time step size is $h = 0.001 \text{s}$ and the simulated duration is $T = 25 \text{s}$. 
	\begin{table}[htbp]
		\begin{center}
			\begin{tabular}{ | c | c| } 
				\hline
				Hovering& $b_d(t)=\left[ 0, \ 0, \ -3\right]\Tp \textup{(m)}$  \\ 
				Slow Swing & $b_d(t)=\left[ 10 \ \textup{sin}(0.1\pi t), \ 0, \ -3\right]\Tp \textup{(m)}$ \\ 
				Fast Swing& $b_d(t)=\left[ 5 \ \textup{sin}(0.5\pi t), \ 0, \ -3\right]\Tp \textup{(m)}$ \\ 
				High Pitch& $b_d(t)=\left[ 10 \ \textup{sin}(0.5\pi t), \ 10 \ \textup{cos}(0.5\pi t), \ -3\right]\Tp \textup{(m)}$ \\
				\hline
			\end{tabular}
		\end{center}
		\caption{Flight trajectories to be tracked for the comparisons between  LESO, FxTSDO and FFTS-ESO}
		\label{table:ESO Comparison Trajectory}
	\end{table}
	\begin{table}[htbp]
		\begin{center}
			\begin{tabular}{ | c| c| c| } 
				\hline
				$b_\textup{N}$ & $b_\textup{N} = b + \mu_b$ & $ \mu_b \sim P_b = 3e^{-8} $   \\ 
				$v_\textup{N}$ & $v_\textup{N} = v + \mu_v$ & $\mu_v \sim P_v = 3e^{-7} $ \\ 
				$R_\textup{N}$ & $R_\textup{N} = R\textup{exp}(\mu_R) $& $\mu_R \sim P_R = 3e^{-8} $ \\ 
				$\Omega_\textup{N}$ & $\Omega_\textup{N} = \Omega + \mu_\Omega $ & $\mu_\Omega \sim P_\Omega = 3e^{-7} $ \\
				\hline
			\end{tabular}
		\end{center}
		\caption{Measurement noise level in power spectral density for the comparisons between LESO, FxTSDO, and FFTS-ESO}
		\label{table:ESO Comparison Noise Level}
	\end{table}
	\begin{figure}[htbp]
		\centering
		\scriptsize
		\subfloat[Hover]{
			\includegraphics[width=\columnwidth]{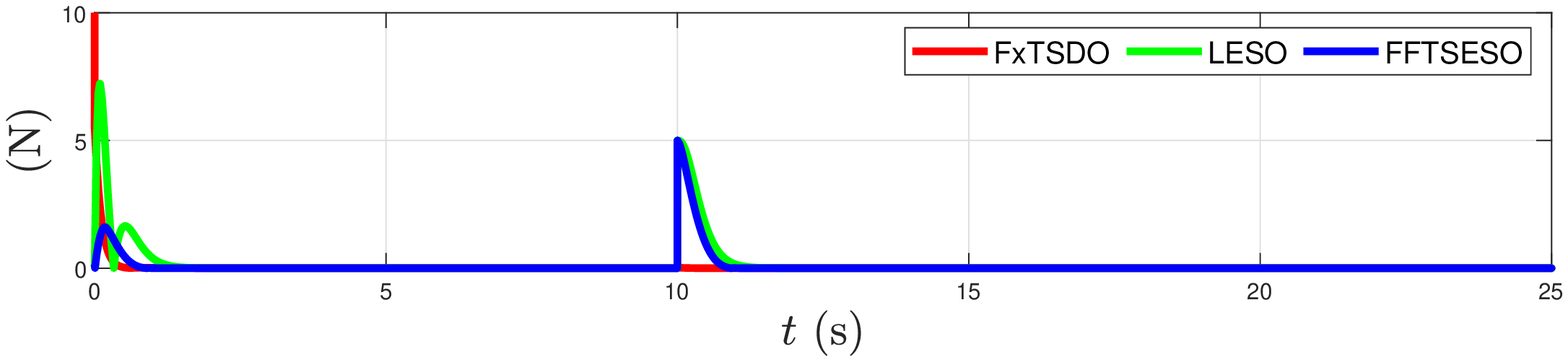}}\\
		\subfloat[Slow swing]{
			\includegraphics[width=\columnwidth]{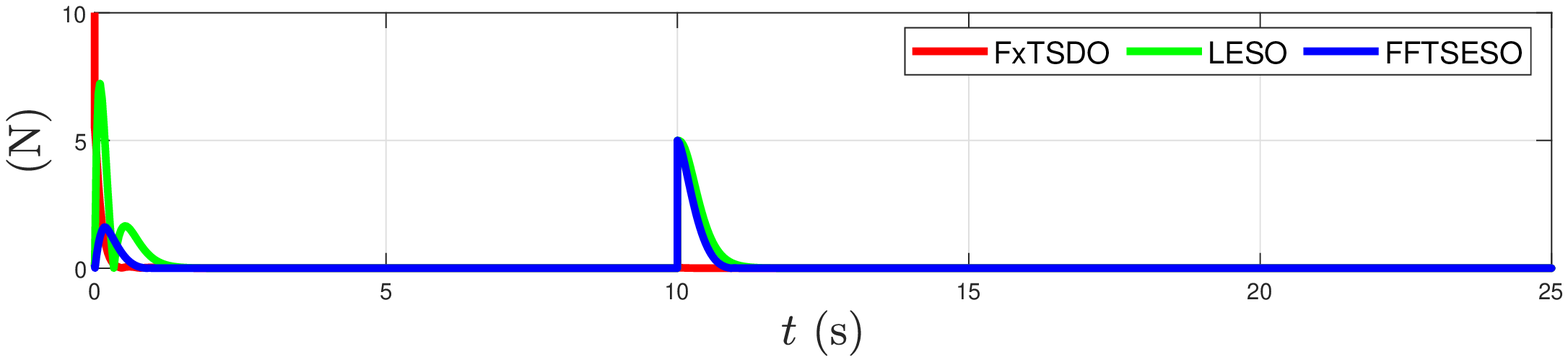}}\\
		\subfloat[Fast swing]{
			\includegraphics[width=\columnwidth]{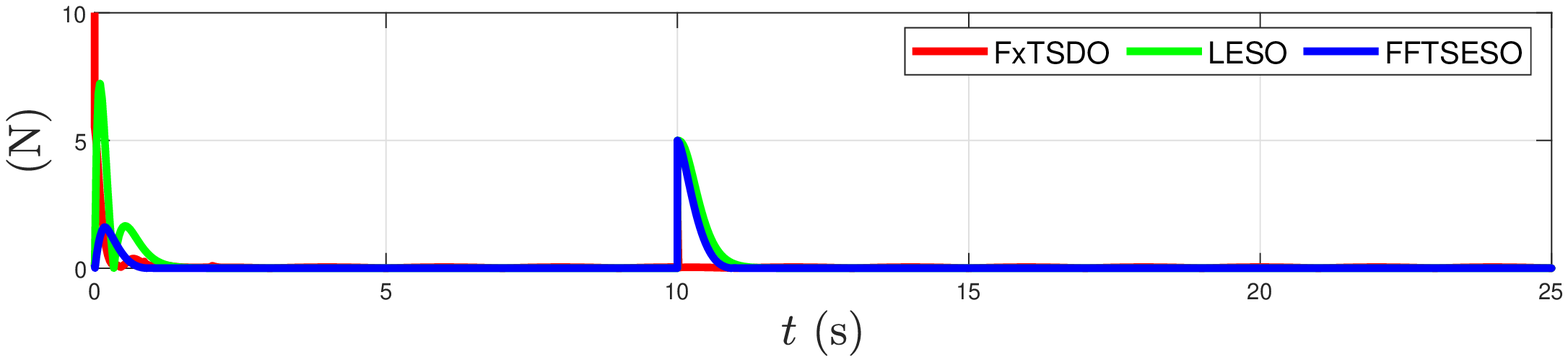}}\\
		\subfloat[High pitch]{
			\includegraphics[width=\columnwidth]{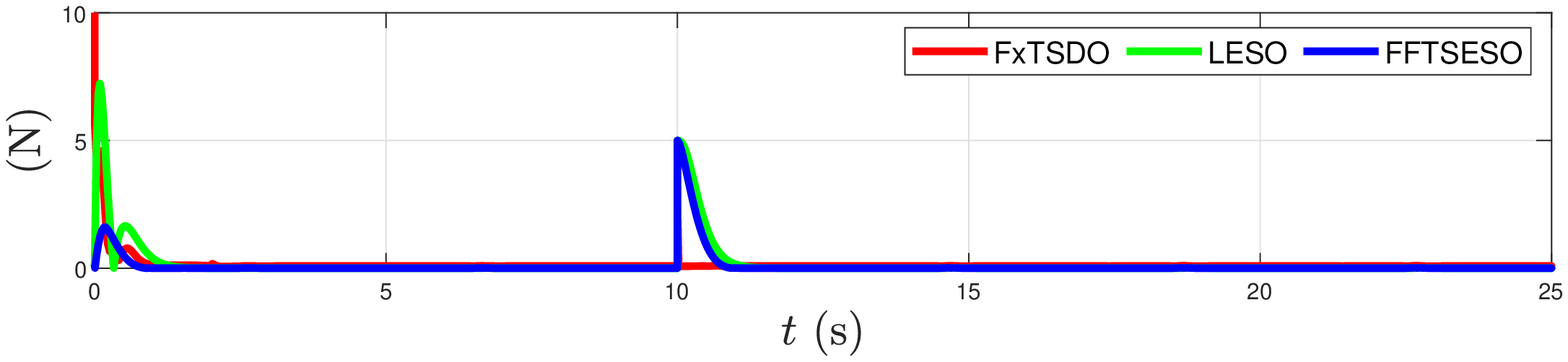}}
		\caption{Disturbance force estimation errors of the multi-rotor UAV from FxTSDO, LESO, and FFTS-ESO, in four different tracking control scenarios without measurement noise. }
		\label{fig:Disturbance Force DOESO Comparison Without Noise}
	\end{figure}
	\begin{figure}[htbp]
		\centering
		\scriptsize
		\subfloat[Hover]{
			\includegraphics[width=\columnwidth]{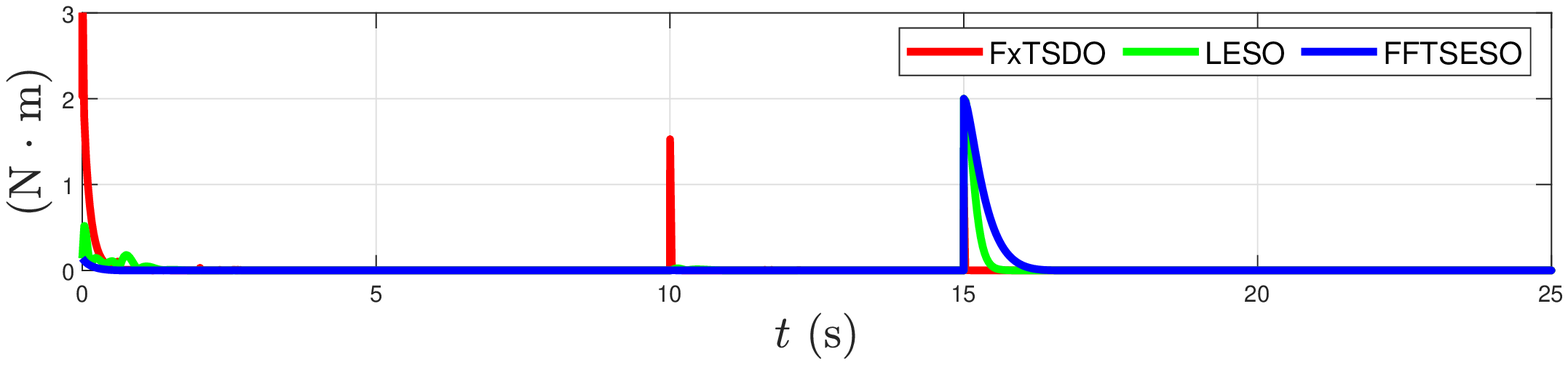}}\\
		\subfloat[Slow swing]{
			\includegraphics[width=\columnwidth]{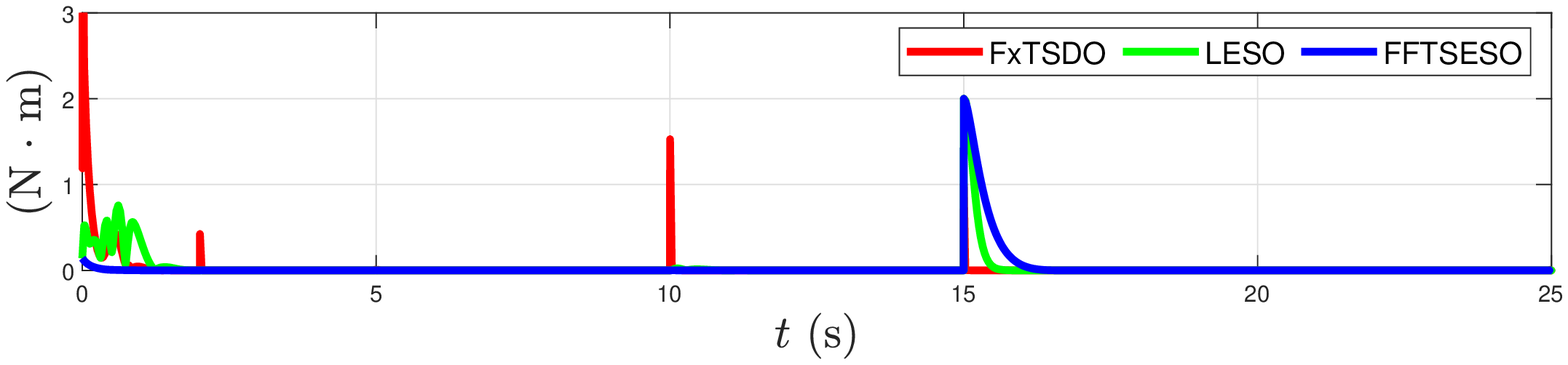}}\\
		\subfloat[Fast swing]{
			\includegraphics[width=\columnwidth]{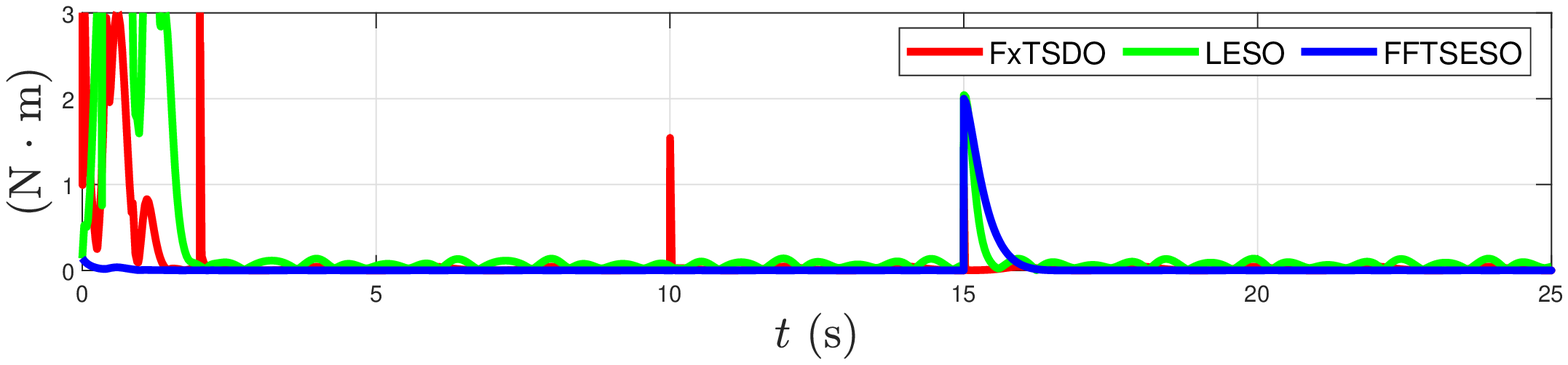}}\\
		\subfloat[High pitch]{
			\includegraphics[width=\columnwidth]{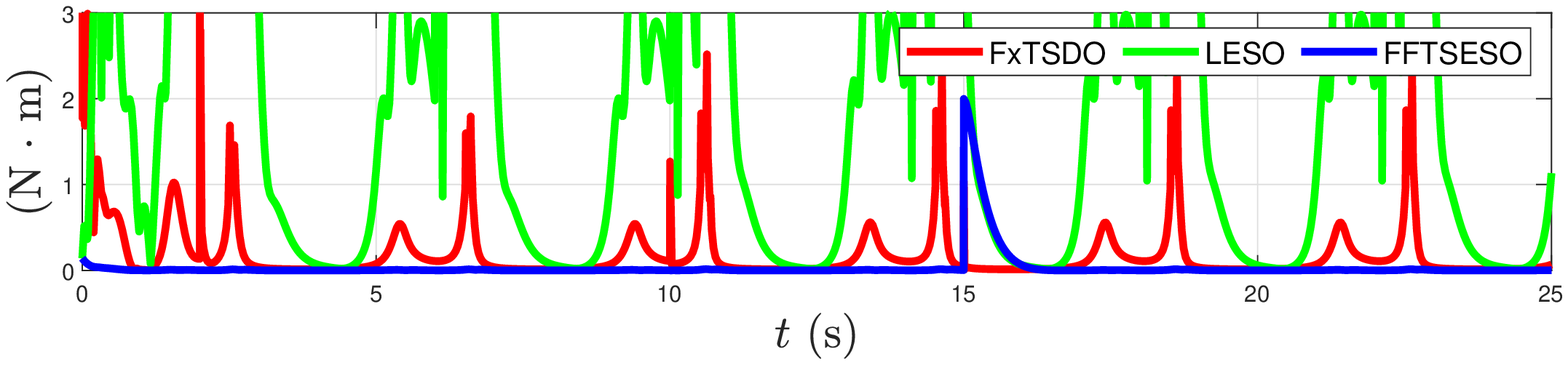}}
		\caption{Disturbance torque estimation errors of the multi-rotor UAV from FxTSDO, LESO, and FFTS-ESO, in four different tracking control scenarios without measurement noise.}
		\label{fig:Disturbance Torque DOESO Comparison Without Noise}
	\end{figure}
	
	The four flight scenarios are the four desired trajectories listed in Table \ref{table:ESO Comparison Trajectory}. `Hovering' is the simplest flight scenario where the aircraft is ordered to hover at a fixed position during the simulation. `High Pitch' is the most complex flight scenario where the aircraft has to pitch up and 
	track a circular trajectory. Since the norm of centripetal acceleration in the `High Pitch' scenario is more than a $g$, the aircraft has to flip over to track the desired trajectory. This desired trajectory with high centripetal acceleration forces the aircraft to go past the 90$^\circ$ pitch 
	singularity of an Euler angle attitude representation. The measurement noise levels are as listed in Table \ref{table:ESO Comparison Noise Level} in terms of power spectral density (PSD). In this set of numerical simulations, the 
	trajectory is tracked by the tracking control system placed in Section \ref{sec:Tracking} without the disturbance rejection term  $\hat{\varphi}_D$ and $\hat{\tau}_D$. The disturbance force and 
	torque in all of the four scenarios in this set of simulations are identical and they are the following step functions: 
	\begin{align*}
		\varphi_{D} (t) &=
		\begin{cases}
			[5, \ 2, \ 0]\T \ \textup{N} & t< 10 \ \textup{s}\\
			[9, \ 5, \ 0]\T \ \textup{N} & t \geq 10 \ \textup{s} \\
		\end{cases}, \\
		\tau_{D} (t) &=
		\begin{cases}
			[2, \ 0, \ 1]\T \ \textup{N}\cdot \textup{m} & t< 15 \ \textup{s} \\
			[4, \ 0, \ 1]\T \ \textup{N}\cdot \textup{m} & t \geq 15 \ \textup{s} \\
		\end{cases}
	\end{align*} 
	
	The parameters for FFTS-ESO in these simulations are $p=1.2, k_{t1}=3, k_{t2}=2,  k_{t3}=2, \kappa_t =0.1$, $k_{a1}=5, k_{a2}=4, k_{a3}=2, \kappa_a =0.3$. The parameters for the tracking control scheme in the simulations are $p=1.2, k_{TP}=5, k_{TD}=16, L_P = I, \kappa_T =2$, $k_{AP}=12, k_{AD}=6, k_{AI}=2, \kappa_A =2, L_A = I$. The gains for FxTSDO and LESO are as given in Liu et al.\cite{liu2022fixed} and Shao et al.\cite{shao2018robust}. 
	In the simulated flight, the initial states of the UAV for all four scenarios are: $R(0) = I, \ \Omega(0) = \left[0, \ 0, \ 0\right]\T \text{rad/s},  b(0) = \left[ 0.01, \ 0, \ 0 \right]\T \text{m}, \ v(0) = \left[5\pi, \ 0, \ 0\right]\T \text{m/s}.$ The initial conditions of the FxTSDO, LESO, and FFTS-ESO, are identical to the pose, velocities and disturbance of the UAV at the initial time in the simulation. 
	
\begin{figure}[htbp]
		\centering
		\scriptsize
		\subfloat[Hover]{
			\includegraphics[width=\columnwidth]{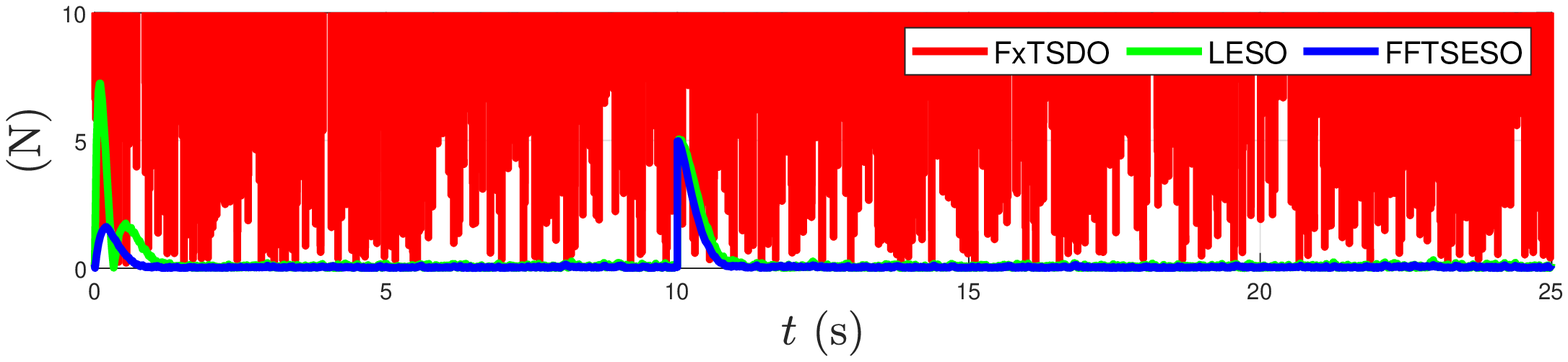}}\\
		\subfloat[Slow swing]{
			\includegraphics[width=\columnwidth]{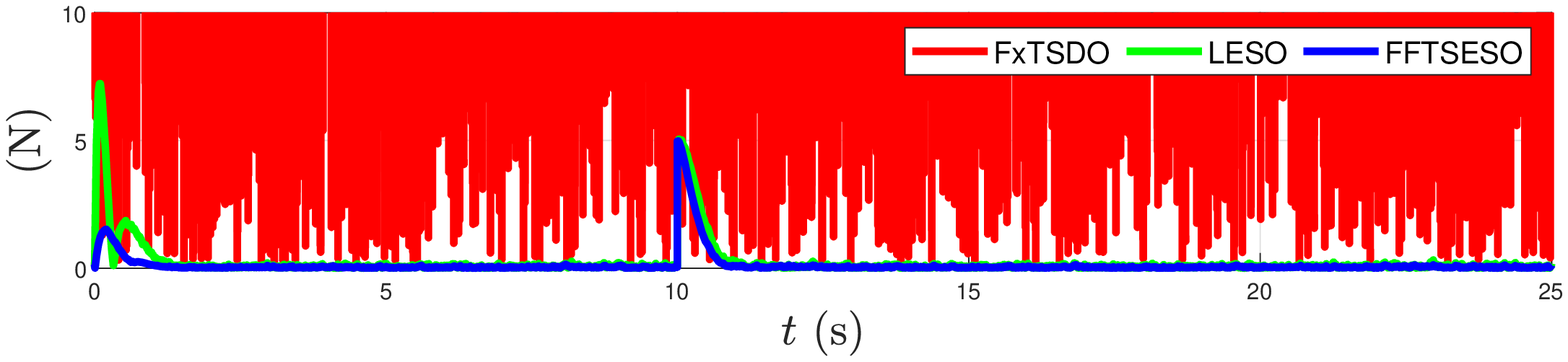}}\\
		\subfloat[Fast swing]{
			\includegraphics[width=\columnwidth]{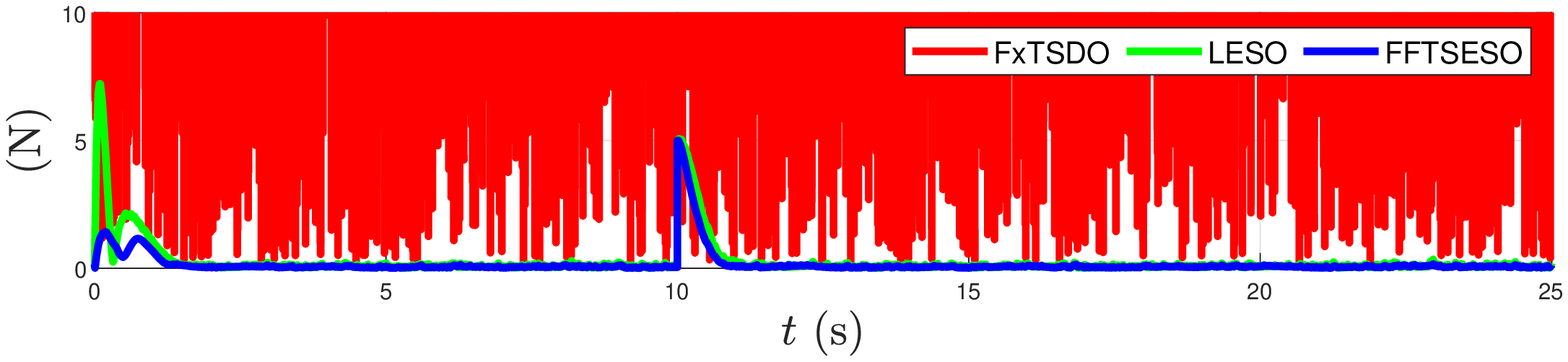}}\\
		\subfloat[High pitch]{
			\includegraphics[width=\columnwidth]{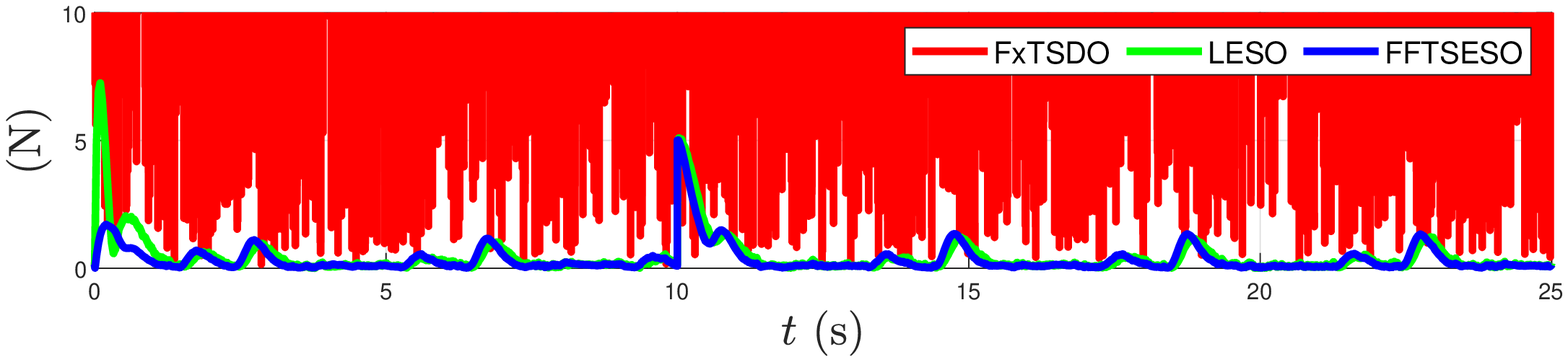}}
		\caption{Disturbance force estimation errors of the multi-rotor UAV from FxTSDO, LESO, and FFTS-ESO, in four different tracking control scenarios with measurement noise.}
		\label{fig:Disturbance Force DOESO Comparison With Noise}
	\end{figure}

\begin{figure}[htbp]
		\centering
		\scriptsize
		\subfloat[Hover]{
			\includegraphics[width=\columnwidth]{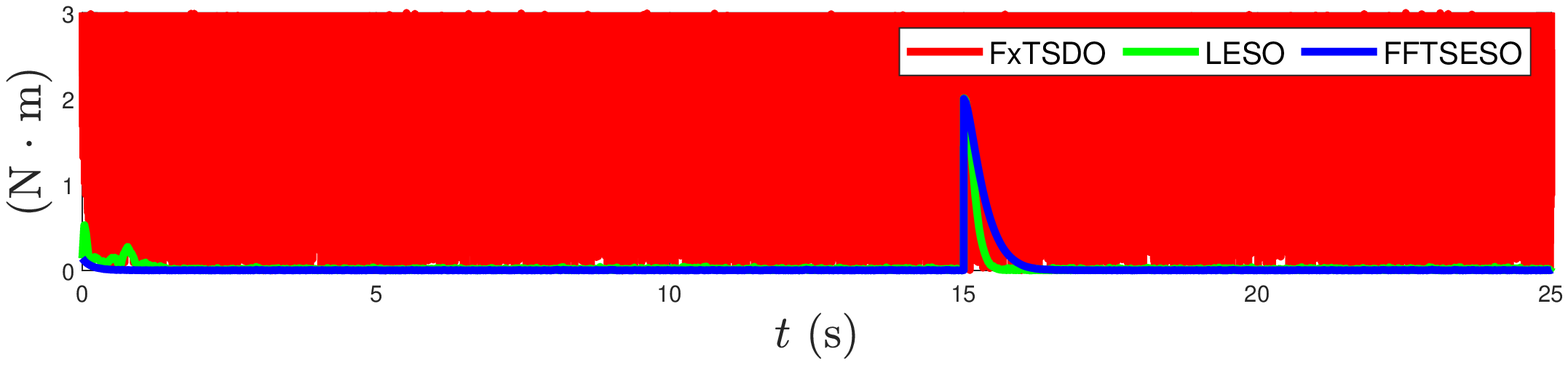}}\\
		\subfloat[Slow swing]{
			\includegraphics[width=\columnwidth]{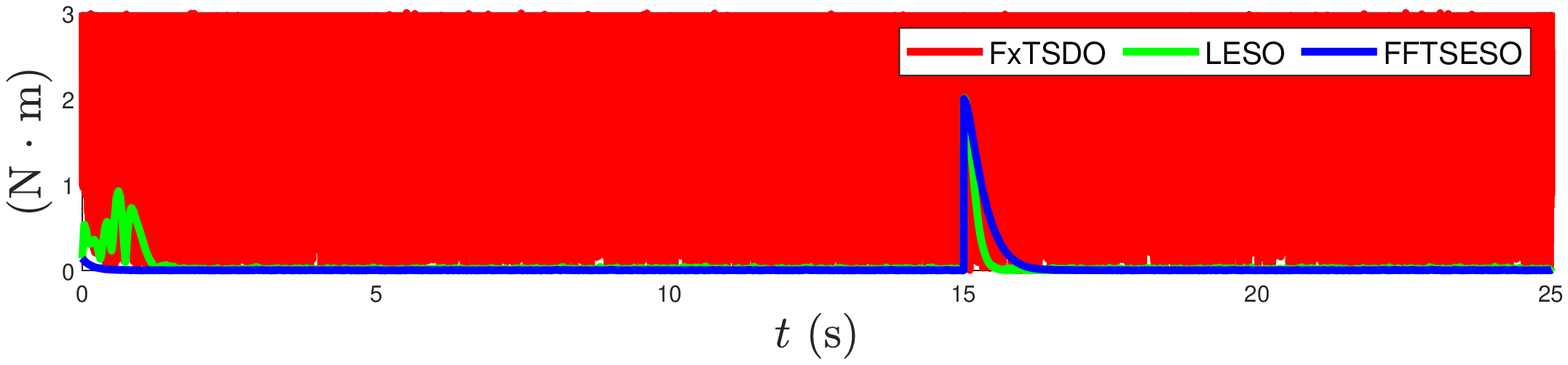}}\\
		\subfloat[Fast swing]{
			\includegraphics[width=\columnwidth]{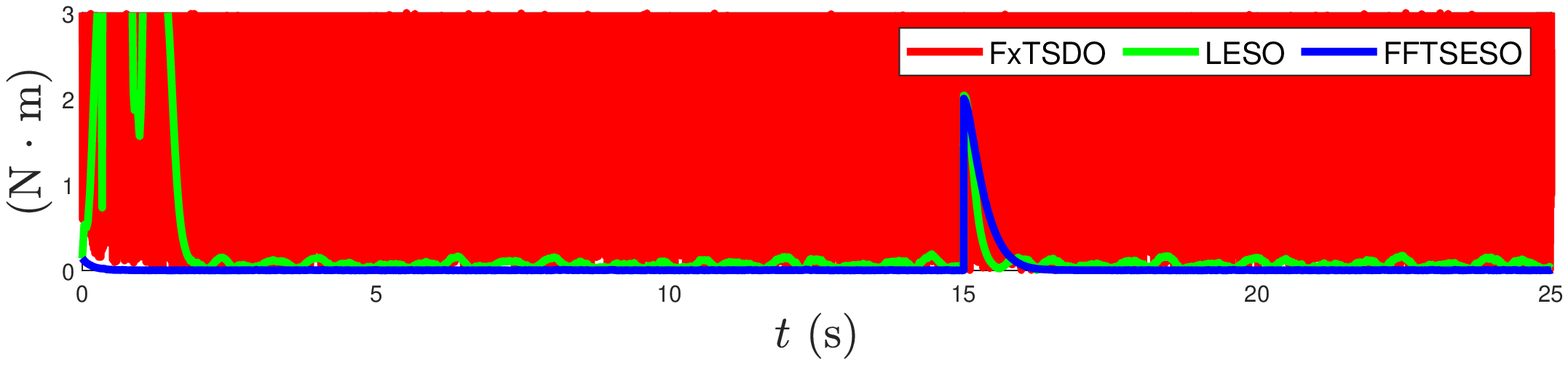}}\\
		\subfloat[High pitch]{
			\includegraphics[width=\columnwidth]{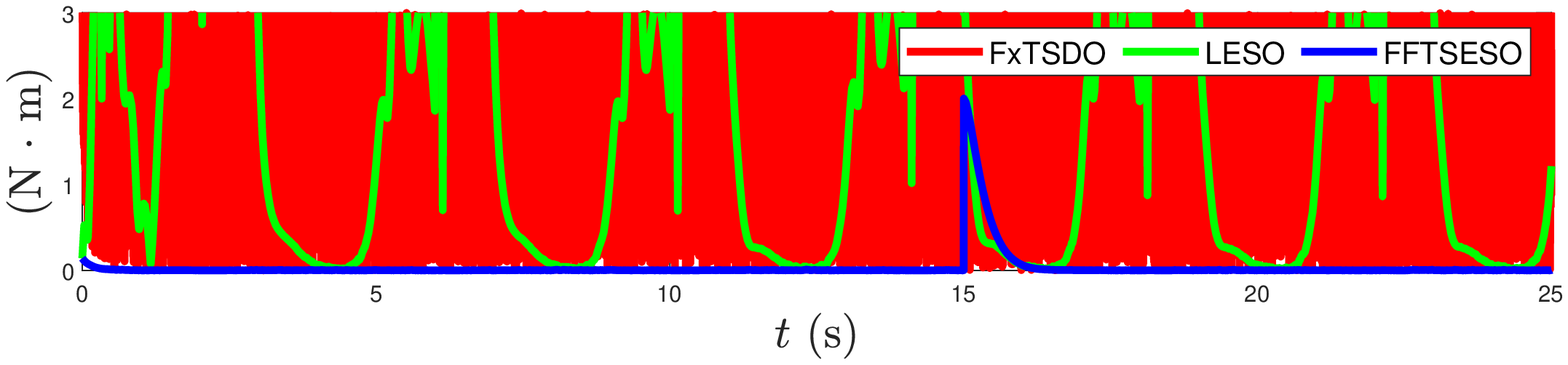}}
		\caption{Disturbance torque estimation errors of the multi-rotor UAV from FxTSDO, LESO, and FFTS-ESO, in four different tracking control scenarios with measurement noise.}
		\label{fig:Disturbance Torque DOESO Comparison With Noise}
	\end{figure}

	We present the simulation results in four sets of figures. Fig.  \ref{fig:Disturbance Force DOESO Comparison Without Noise} and \ref{fig:Disturbance Torque DOESO Comparison Without Noise} present the disturbance force and torque estimation errors respectively, from FxTSDO, LESO and FFTS-ESO in the flight scenarios described in Table \ref{table:ESO Comparison Trajectory} with noise-free measurements. Fig. \ref{fig:Disturbance Force DOESO Comparison With Noise} and \ref{fig:Disturbance Torque DOESO Comparison With Noise} present the disturbance estimation errors from these schemes for the flight trajectories in Table \ref{table:ESO Comparison Trajectory}, in the presence of  measurement noise levels as described in Table \ref{table:ESO Comparison Noise Level}.
	
	Fig. \ref{fig:Disturbance Force DOESO Comparison Without Noise} shows the disturbance force estimation errors from the three schemes with noise-free measurements. Although the disturbance force estimation error from FxTSDO shows significant initial transient, the results from Fig. \ref{fig:Disturbance Force DOESO Comparison Without Noise} indicates that with noise-free measurement, the disturbance force estimations from these three schemes converge to the origin in all four flight scenarios. The transients at $t=15$ s are from the step-function disturbance force $\varphi_D$, whose step time is $t=15$ s.
	Fig. \ref{fig:Disturbance Torque DOESO Comparison Without Noise} shows the disturbance torque estimation errors from the three schemes with noise-free measurement. In Fig. \ref{fig:Disturbance Torque DOESO Comparison Without Noise}, we observe that when $t=10 \ \text{s}$, high transients appears in the disturbance torque estimation error from FxTSDO. 
	
	Despite the initial transients, the disturbance torque estimation errors from all three schemes converge to the origin in 'Hovering' and 'Slow swing' scenarios. However, in 'Fast swing' and 'High pitch' scenarios, the disturbance torque estimation errors from LESO and FxTSDO diverge. As is stated in Section \ref{sec:Intro}, since the LESO uses Euler-angle to represent attitude for disturbance torque estimation, it experiences a singularity in attitude representation when the UAV tracks the 'Fast swing' and 'High Pitch' trajectories. Thus, in these two scenarios, the singularity in the attitude representation destabilizes the disturbance torque estimation error of LESO. 
	Fig. \ref{fig:Disturbance Force DOESO Comparison With Noise} and \ref{fig:Disturbance Torque DOESO Comparison With Noise} present the disturbance force and disturbance torque estimation errors respectively, from the three schemes with identical noisy measurements as given in Table \ref{table:ESO Comparison Noise Level}. As is stated in Remark \ref{rem:Compare}, we observe that with measurement noise, FxTSDO is not capable of providing any meaningful disturbance estimation.  In 'Fast swing' and 'High pitch' scenarios, the disturbance torque estimation errors from LESO diverge from the origin. 
	
	To summarize, Fig. \ref{fig:Disturbance Force DOESO Comparison Without Noise}, \ref{fig:Disturbance Torque DOESO Comparison Without Noise}, \ref{fig:Disturbance Force DOESO Comparison With Noise}, and \ref{fig:Disturbance Torque DOESO Comparison With Noise} show that the FFTS-ESO has satisfactory disturbance estimation performance and outperforms the LESO and FxTSDO when the UAV experiences large pose changes and has noisy measurements.
	
	\subsection{Simplified simulation with LGVI}\label{subsec:LGVI}
	In this subsection, we present a set of numerical simulation results on the FFTS-ADRC scheme, based on the geometric integrator LGVI\cite{nordkvist2010lie}. The inertia and mass of the aircraft is $  J=\textup{diag}([0.0820,0.0845,0.1377]) \ \textup{kg}\cdot \textup{m}^2, \quad m = 4.34 \ \textup{kg}. $ The motion of the UAV is numerically integrated in discrete time with an LGVI. The simulated time duration is $T=5 \textup{s}$, with a time step size of $\Delta t=0.005\textup{s}$. The ESO and control gains of the implemented FFTS-ADRC are: $k_{t1}=5;  k_{t2}=5; k_{t3}=3; \kappa_t=2; k_{TD} = 4; \ L_T=I; k_{TP} = 2; \kappa_T = 2, k_{a1}=5;  k_{a2}=6; k_{a3}=3; \kappa_a=1.5; k_{AD} = 3; \ L_A=0.5I; k_{AP} = 3; k_{AI} = 0.1;  \kappa_A = 2; p=1.2.$
	The desired trajectory of the simulation is given by: $b_d(t)=\left[ 2 \ \textup{sin}(\pi t), \ 2 t, \ 2\textup{sin}(\pi t) \right]\T \textup{m},$
	resulting in the aircraft experiencing a singularity in the Euler angle attitude representation.
	
	The dynamic disturbance force and torque in all of the four scenarios in this set of simulations are identical and they are the following functions: 
	\begin{align*}
		\varphi_{D} (t) &=
		\begin{bmatrix}
			&50+6\textup{sin}\left(\frac{\pi t}{2}\right)+\frac{1}{2}\textup{sin}\left(\pi t\right) \\
			&50+3\textup{sin}\left(\frac{\pi t}{2}\right)+\frac{1}{5}\textup{sin}\left(\pi t\right) \\
			&20 
		\end{bmatrix} \ \textup{N} \\
		\tau_{D} (t) &= 
		\begin{bmatrix}
			&5+\frac{1}{2}\textup{sin}\left(\frac{\pi t}{2}\right)+\frac{1}{10}\textup{sin}\left(\pi t\right)\\
			&3+\textup{sin}\left(\frac{\pi t}{2}\right)+\frac{1}{20}\textup{sin}\left(\pi t\right)\\
			&-3
		\end{bmatrix} \ \textup{N}\cdot \textup{m}  
	\end{align*} 
	
	In the simulated flight, the initial states of the UAV for all four scenarios are: $R(0) = I, \ \Omega(0) = \left[0, \ 0, \ 0\right]\T \text{rad/s}, b(0) = \left[ 0, \ 0, \ 3 \right]\T \text{m}, \ v(0) = \left[2\pi, \ 0, \ 0\right]\T \text{m/s}. $
	The initial conditions of the FFTS-ESO are identical to the pose, velocities and disturbance inputs of the UAV at the initial time for this simulation. 
	
\begin{figure}[htbp]
		\centering
		\scriptsize
		\subfloat[Without rejection]
		{\includegraphics[width=0.4\columnwidth]{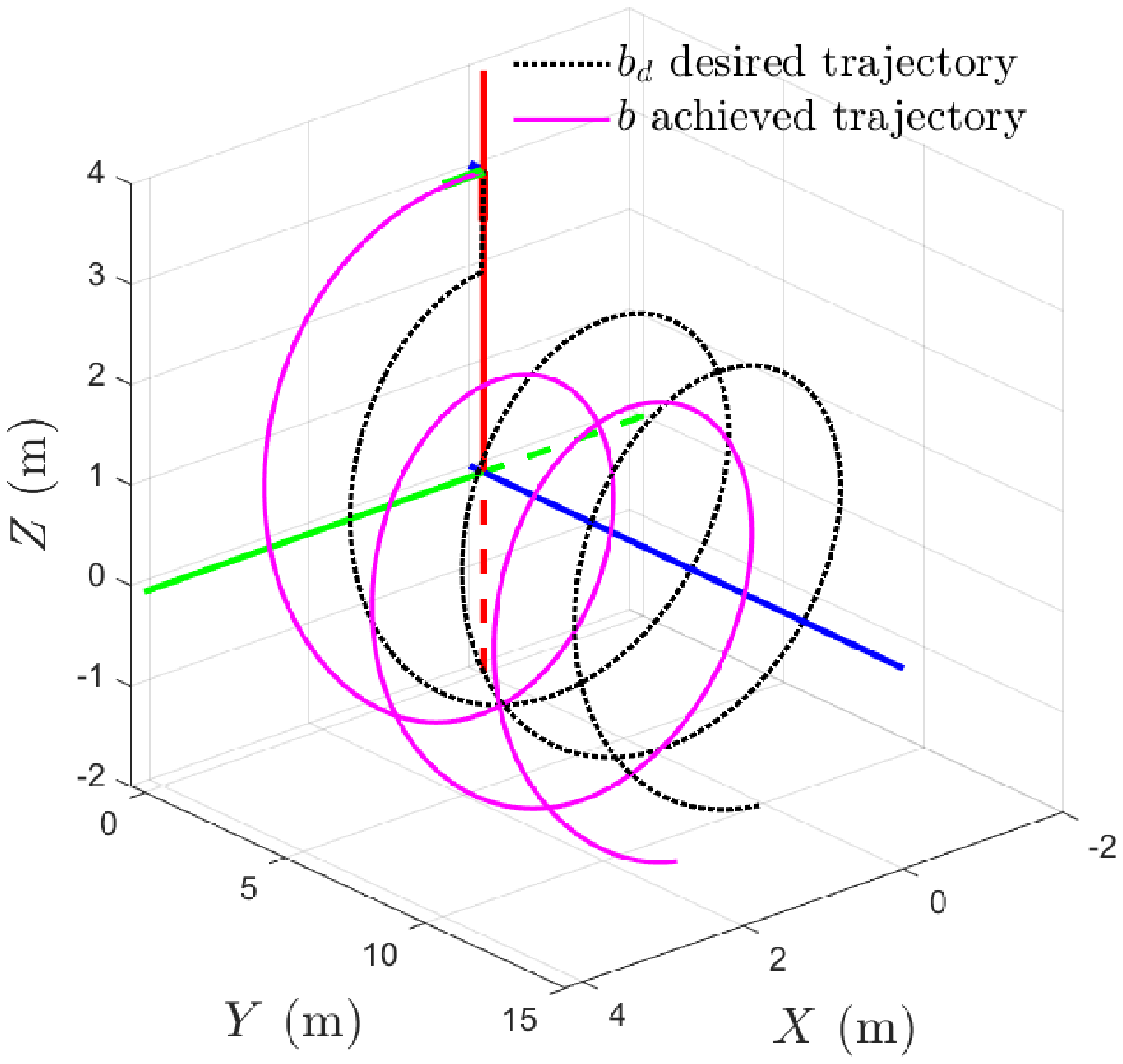}}
		\subfloat[Disturbance force rejection]
		{\includegraphics[width=0.4\columnwidth]{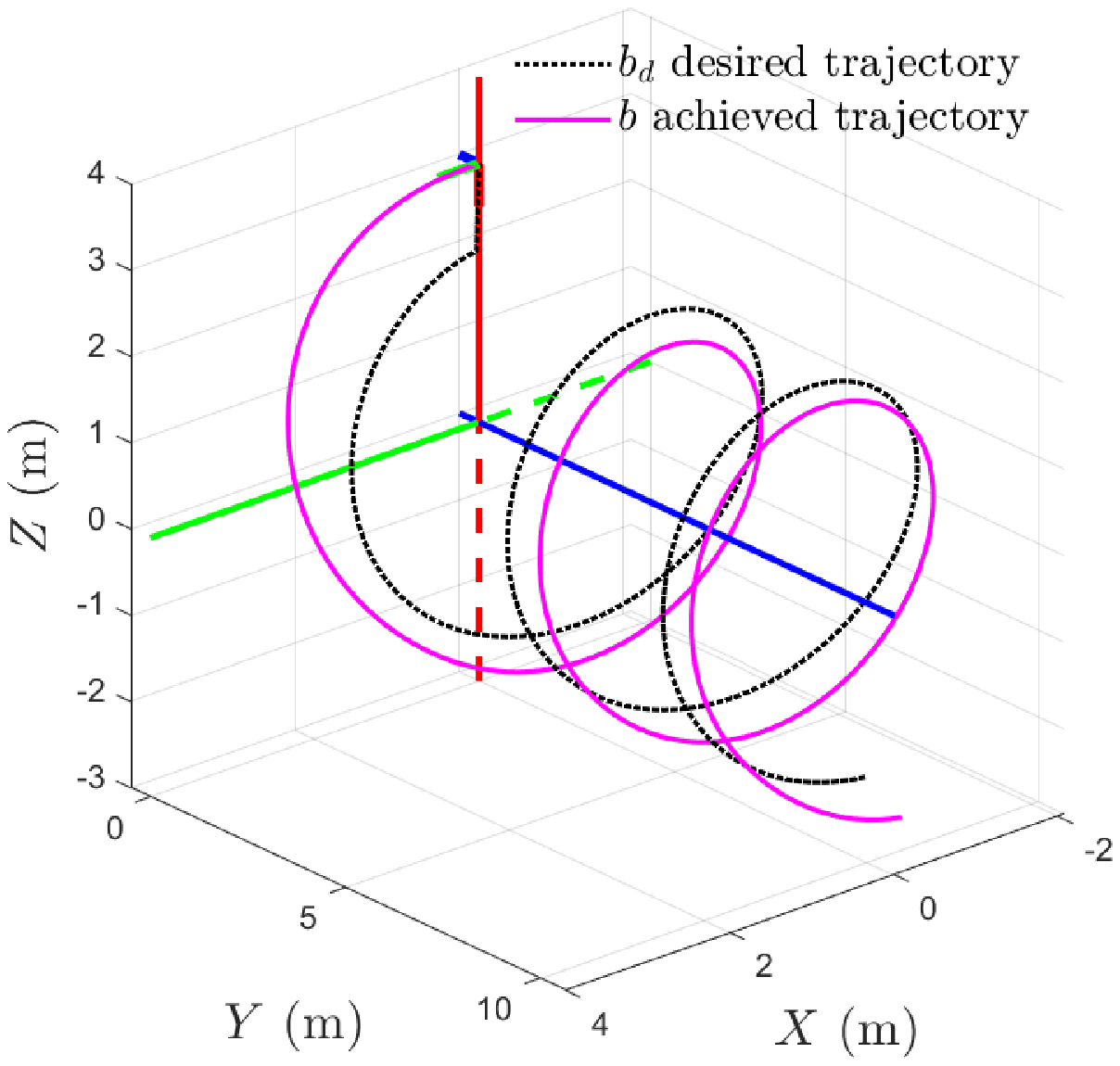}}\\
		\subfloat[Disturbance torque rejection]
		{\includegraphics[width=0.4\columnwidth]{FFTS-ADRC/ADRCComparison_TransADRC.eps}}
		\subfloat[With both rejection]
		{\includegraphics[width=0.4\columnwidth]{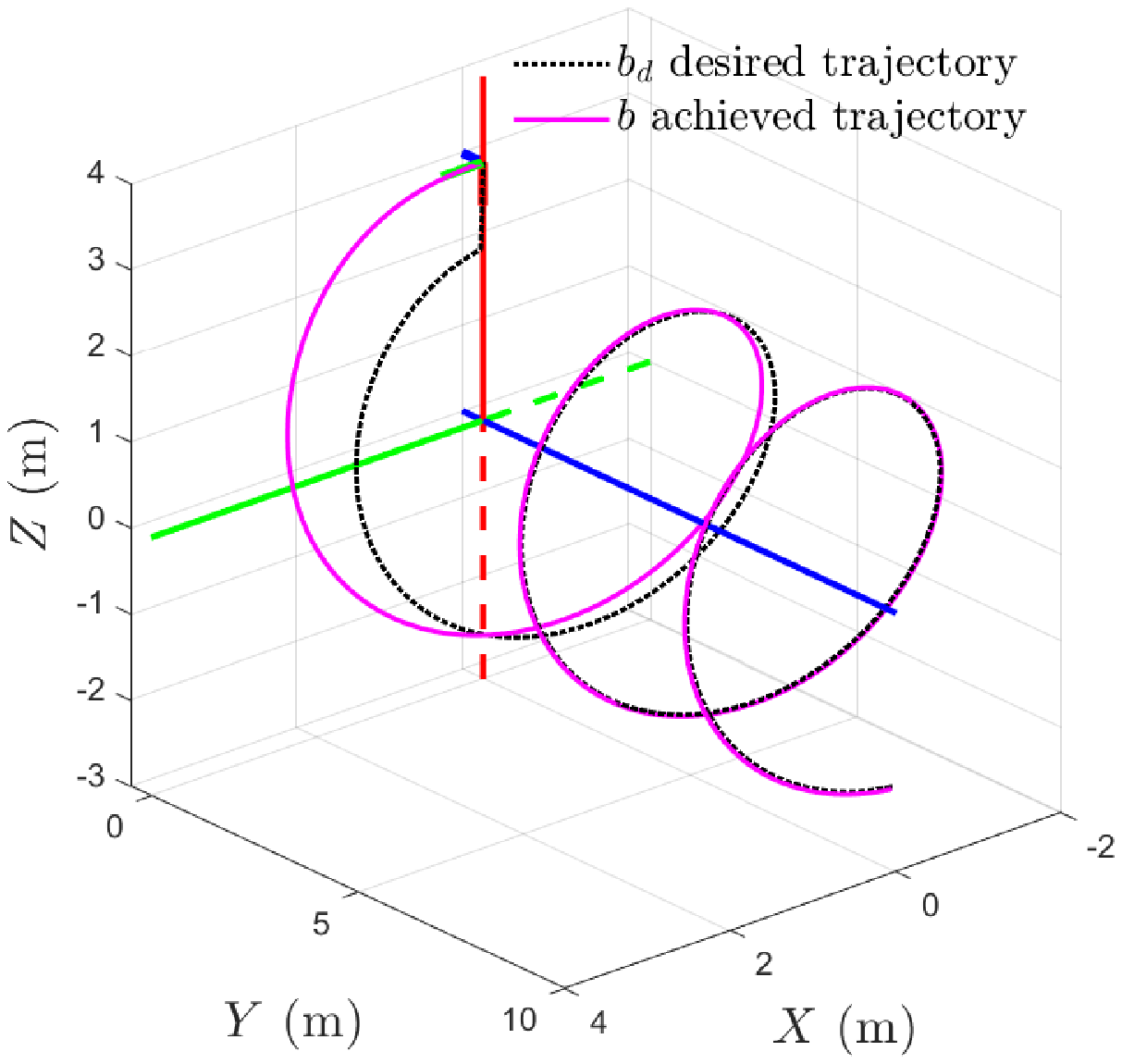}}
		\caption{The tracked trajectories of the proposed FFTS-ADRC with different configurations of disturbance rejection}
		\label{fig:FFTS-ADRC LGVI Traj}
	\end{figure}

\begin{figure}[htbp]
		\centering
		\scriptsize
		\subfloat[Position tracking error]
		{\includegraphics[width=.8\columnwidth]{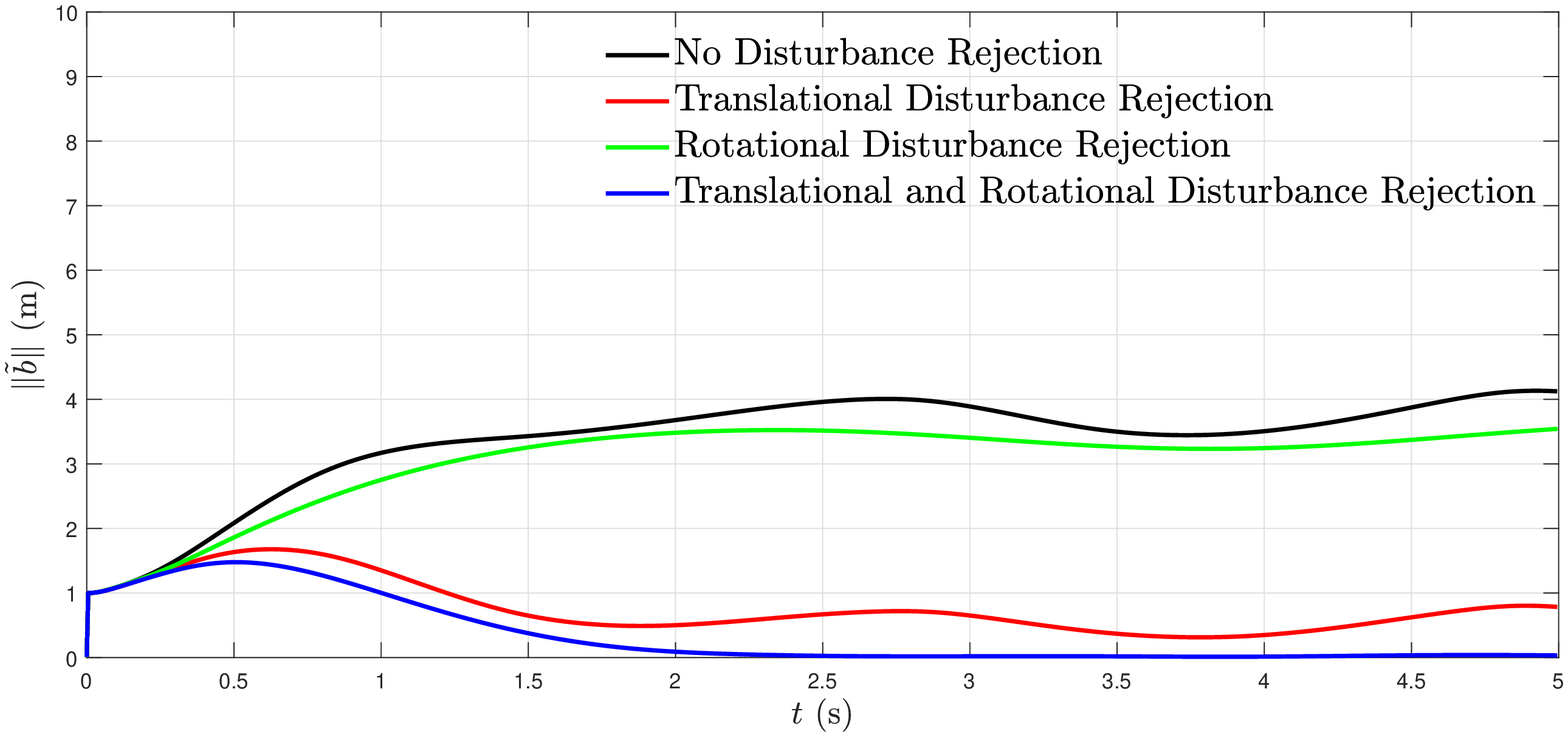}}\\
		\subfloat[Attitude tracking error]
		{\includegraphics[width=.8\columnwidth]{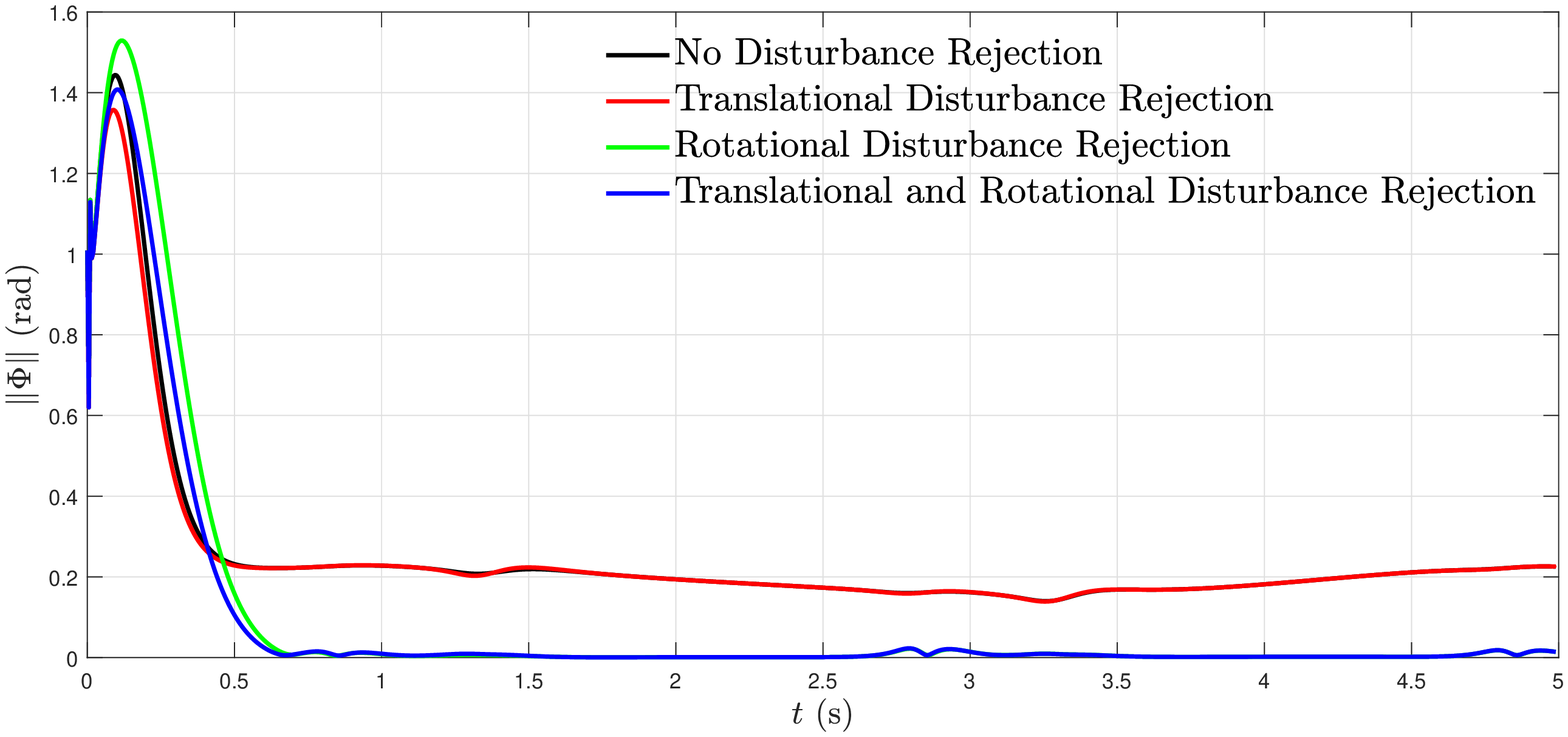}}
		\caption{Position and attitude tracking error of the FFTS-ADRC with different configurations of disturbance rejection}
		\label{fig:FFTS-ADRC Error}
	\end{figure}

	Four sets of simulation results, namely simulation without disturbance rejection, with only disturbance force rejection, with only disturbance torque rejection, and with both disturbance force and disturbance torque rejection, are included in this section to validate the control performance of the proposed FFTS-ADRC scheme. The results are presented in Fig. \ref{fig:FFTS-ADRC LGVI Traj} and Fig. \ref{fig:FFTS-ADRC Error}.
	From Fig. \ref{fig:FFTS-ADRC LGVI Traj}, we observe that all of the trajectories of the simulated flights converge to a neighborhood near the desired trajectory. Among these trajectories of the simulated flights, the one with both disturbance force and disturbance torque rejection is the one closest to the desired trajectory.
	
	Fig. \ref{fig:FFTS-ADRC Error} shows the results for the attitude and position tracking errors. The attitude tracking error $Q$ is parameterized by the principal rotation angle  
	$\Phi = \textup{cos}^{-1} \left(\frac{1}{2}(\textup{tr}(Q)-1)\right).$
	The position tracking error is defined as the norm of $\tilde{b}$. In Fig. \ref{fig:FFTS-ADRC Error}, the $\|\Phi\|$ curve of the simulated flight without disturbance rejection and the one with only disturbance force rejection are almost identical to each other. Fig. \ref{fig:FFTS-ADRC Error} indicates that the simulated flight with both disturbance force and disturbance torque rejection has the best control performance.
	
	\section{Conclusion}\label{sec:Conclusion}
	In this article, a feedback tracking geometric control scheme using a FFTS-ESO for disturbance rejection is designed for rotorcraft UAV that have a body-fixed thrust direction and three-axis attitude control. The resulting FFTS-ADRC scheme enables such rotorcraft UAV to perform large maneuvers in the presence of aerodynamic uncertainties. The UAV is modeled as an underactuated system on the tangent bundle of the six-dimensional Lie group of rigid body motions, $\SE$. The proposed ESO scheme is developed based on the HC-FFTSD, which is similar to the STA used in sliding mode designs, to obtain fast finite-time stability with higher tunability of the settling time compared to other FTS schemes. The ADRC scheme on $\SE$, which utilizes the estimated disturbances from the ESO, is then incorporated to achieve FFTS tracking errors under constant disturbances and ultimate boundedness of tracking errors for time-varying disturbances. The Lyapunov stability analysis presented in this article for both ESO scheme and tracking control scheme proves the finite-time stability and robustness of the overall ADRC on $\SE$. Two sets of numerical simulations are conducted.  The first set of numerical simulation results present the stable performance of the FFTS-ESO scheme in estimating external force and torque disturbances acting on the UAV in different scenarios. The behavior of the FFTS-ESO is compared with two state-of-the-art observers for disturbance estimation. Using a realistic set of data for several simulated flight scenarios of a rotorcraft UAV, numerical simulations show that the FFTS-ESO, unlike the LESO and FxTSDO, is always stable and its convergence is robust to measurement noise and pose singularities. The proposed FFTS-ADRC scheme is numerically implemented by a geometric integrator for a rotorcraft UAV model and numerical simulations are carried out to validate the developed FFTS-ESO and FFTS-ADRC schemes. The numerical results also demonstrate the stable performance of the FFTS-ADRC when the UAV carries out large maneuvers that lead to kinematic singularities in Euler angle attitude representation. 
	\appendix
	\begin{Proof*}
		{\em
			\textup{Represent $x$ as a linear combination of $\mu$ and $\nu$: }
			\begin{align}\label{eqn:x representation}
				x = c_1 \mu +c_2\nu,  
			\end{align}
			\textup{where $\nu$ is a vector perpendicular to $\mu$, such that $\mu\T\nu=0$. Next, define two non-zero scalars, $c_1, c_2$. Using \eqref{eqn:x representation}, express $Y$ in Lemma \ref{lem:Inequality Noise Robustness} in coordinates $(c_1,c_2)$: }
			\begin{align*}
				Y = \frac{c_1\mu+c_2\nu}{\left(c_1^2\|\mu\|^2 +c_2^2\|\nu\|^2\right)^\alpha} - \frac{(1+c_1)\mu+c_2\nu}{\left[(1+c_1)^2\|\mu\|^2 +c_2^2\|\nu\|^2\right]^\alpha}.
			\end{align*}
			\textup{Thereafter, we obtain its partial derivatives with respect to these coordinates:}
			\begin{align}\label{eqn:Y partial c1}
				\begin{split}
					\frac{\partial Y}{\partial c_1} &= \frac{\mu}{\left(c_1^2\|\mu\|^2 +c_2^2\|\nu\|^2\right)^\alpha}\\
					&- \frac{2\alpha c_1\|\mu\|^2(c_1\mu+c_2\nu)}{\left(c_1^2\|\mu\|^2 +c_2^2\|\nu\|^2\right)^{\alpha+1}}  \\
					&- \frac{\mu}{\left[(1+c_1)^2\|\mu\|^2 +c_2^2\|\nu\|^2\right]^\alpha} \\
					&+ \frac{2\alpha( 1+c_1)\|\mu\|^2\left[(1+c_1)\mu+c_2\nu\right]}{\left[(1+c_1)^2\|\mu\|^2 +c_2^2\|\nu\|^2\right]^{\alpha+1}},
				\end{split}
			\end{align}
			\begin{align}\label{eqn:Y partial c2}
				\begin{split}
					\frac{\partial Y}{\partial c_2} &= \frac{\nu}{\left(c_1^2\|\mu\|^2 +c_2^2\|\nu\|^2\right)^\alpha} \\
					&- \frac{2\alpha c_2\|\nu\|^2(c_1\mu+c_2\nu)}{\left(c_1^2\|\mu\|^2 +c_2^2\|\nu\|^2\right)^{\alpha+1}} \\
					&- \frac{\nu}{\left[(1+c_1)^2\|\mu\|^2 +c_2^2\|\nu\|^2\right]^\alpha} \\
					&+ \frac{2\alpha c_2\|\nu\|^2\left[(1+c_1)\mu+c_2\nu\right]}{\left[(1+c_1)^2\|\mu\|^2 +c_2^2\|\nu\|^2\right]^{\alpha+1}}.
				\end{split}
			\end{align}
			\textup{Thereafter, we employ the fact that the local maxima of $Y\Tp Y$ satisfy:}
			\begin{align*}
				\begin{split}
					\frac{\partial }{\partial c_1} (Y\Tp Y) = \frac{\partial }{\partial c_2} (Y\Tp Y) = 0, 
				\end{split}
			\end{align*}
			we obtain the following equivalent conditions for the maxima:
			\begin{align}
				&\nu\Tp \frac{\partial Y}{\partial c_1} = \mu\Tp \frac{\partial Y}{\partial c_2} = 0, 	\label{eqn:Y partial equation 01}\\
				&\mu\Tp \frac{\partial Y}{\partial c_1} = 0, 											\label{eqn:Y partial equation 02}\\
				&\nu\Tp \frac{\partial Y}{\partial c_2} = 0. 											\label{eqn:Y partial equation 03}
			\end{align}
			Substituting \eqref{eqn:Y partial c1} and \eqref{eqn:Y partial c2} into \eqref{eqn:Y partial equation 01}, we obtain:
			\begin{align}\label{eqn:Y partial equation 1a} 
				\begin{split}
					&\nu\Tp \frac{\partial Y}{\partial c_1} = \mu\Tp \frac{\partial Y}{\partial c_2} =  0, \\
					&\Longleftrightarrow -\frac{2\alpha c_1 c_2 \|\mu\|^2\|\nu\|^2 }{\left(c_1^2\|\mu\|^2 +c_2^2\|\nu\|^2\right)^{\alpha+1}} \\
					&+ \frac{2\alpha (1+c_1) c_2 \|\mu\|^2\|\nu\|^2}{\left[(1+c_1)^2\|\mu\|^2 +c_2^2\|\nu\|^2\right]^{\alpha+1}} =0, \\
					&\Longrightarrow  \ c_1\left[(1+c_1)^2\|\mu\|^2 +c_2^2 \|\nu\|^2\right]^{\alpha+1} \\
					&= (1+c_1)\left[c_1^2\|\mu\|^2 +c_2^2 \|\nu\|^2\right]^{\alpha+1}, \\
				\end{split}
			\end{align}
			Substituting \eqref{eqn:Y partial c1} and \eqref{eqn:Y partial c2} into \eqref{eqn:Y partial equation 02}, we obtain:
			\begin{align}\label{eqn:Y partial equation 2}
				\begin{split}
					&\mu\Tp \frac{\partial Y}{\partial c_1} = 0, \\
					&\Longrightarrow  \frac{(1-2\alpha \|\mu\|^2 c_1^2)\|\mu\|^2}{\left(c_1^2\|\mu\|^2 +c_2^2\|\nu\|^2\right)^{\alpha+1}} \\
					&- \frac{\left[ 1-2\alpha (1+c_1)^2\|\mu\|^2 \right]\|\mu\|^2}{\left[(1+c_1)^2\|\mu\|^2 +c_2^2\|\nu\|^2\right]^{\alpha+1}} = 0, \\
					&\Longleftrightarrow  (1+c_1)^2 = c_1^2, \Longleftrightarrow  c_1 = -\frac{1}{2}.
				\end{split}
			\end{align}
			Substituting \eqref{eqn:Y partial c1} and \eqref{eqn:Y partial c2} into \eqref{eqn:Y partial equation 03}, we obtain:
			\begin{align}\label{eqn:Y partial equation 3}
				\begin{split}
					&\nu\Tp \frac{\partial Y}{\partial c_2} = 0, \\
					&\Longrightarrow  \frac{(1-2\alpha \|\nu\|^2 c_2^2)\|\nu\|^2}{\left(c_1^2\|\mu\|^2 +c_2^2\|\nu\|^2\right)^{\alpha+1}} \\
					&- \frac{(1-2\alpha \|\nu\|^2 c_2^2)\|\nu\|^2}{\left[(1+c_1)^2\|\mu\|^2 +c_2^2\|\nu\|^2\right]^{\alpha+1}} = 0, \\
					\Longleftrightarrow & (1+c_1)^2 = c_1^2, 
					\Longleftrightarrow  c_1 = -\frac{1}{2}. \\
				\end{split}
			\end{align}
			\eqref{eqn:Y partial equation 1a} does not give a real solution for $\alpha \in ]0,1/2[  $. Thus, we conclude that the only solution to \eqref{eqn:Y partial equation 01}, \eqref{eqn:Y partial equation 02}, \eqref{eqn:Y partial equation 03} is given by $c_1 = -1/2, c_2 = 0$. 
			Thus, the only critical value of $Y\Tp Y$ is obtained when $x=-\mu/2$.
			Further, we conclude that the global maximum of $Y\Tp Y$ is at $x=-\mu/2$ because it is positive definite in $Y$. Therefore, we do not need an analysis of the Hessian matrix of $Y\T Y$ as a function of $(c_1,c_2)$. 
		} 
		\qed
	\end{Proof*}

	\bibliographystyle{IEEEtran}
	\bibliography{ref_TCST.bib}

\end{document}